\documentclass[10pt]{amsart}
\usepackage[pagebackref]{hyperref}

\renewcommand{\d}{\mathrm{d}}%exterior derivative

\newcommand{\ts}{\textstyle }
\newcommand{\ds}{\displaystyle }

\newcommand{\bbR}{{\mathbb R}}
\newcommand{\bbC}{{\mathbb C}}

\newcommand{\bbQ}{{\mathbb Q}}

\newcommand{\bbT}{{\mathbb T}}
\newcommand{\bbZ}{{\mathbb Z}}

\newcommand{\E}{{\mathrm e}}
\newcommand{\iC}{{\mathrm i}}
\newcommand{\GL}{\operatorname{GL}}

\newcommand{\Un}{\operatorname{U}}

\newcommand{\Ric}{\operatorname{Ric}}

\newcommand{\Lie}{\operatorname{\mathsf{L}}}

\DeclareMathOperator{\tr}{tr}

\newcommand{\p}{\partial}

\newcommand{\cI}{\mathcal{I}}

\newcommand{\hs}{\mathsf{h}}

\newcommand{\ks}{\mathsf{k}}

\newcommand{\eut}{\operatorname{\mathfrak{t}}}
\newcommand{\euu}{\operatorname{\mathfrak{u}}}

\newcommand{\bi}{{\bar{\imath}}}
\newcommand{\bj}{{\bar{\jmath}}}
\newcommand{\bk}{{\bar k}}
\newcommand{\bl}{{\bar l}}

\newcommand{\wbar}{{\bar w}}
\newcommand{\zbar}{{\bar z}}

\newcommand{\la}{\langle}
\newcommand{\ra}{\rangle}

\newcommand{\w}{{\mathchoice{\,{\scriptstyle\wedge}\,}%
{{\scriptstyle\wedge}}{{\scriptscriptstyle\wedge}}%
{{\scriptscriptstyle\wedge}}}}
\newcommand{\lhk}{\mathbin{\hbox{\vrule height1.4pt width4pt depth-1pt 
             \vrule height4pt width0.4pt depth-1pt}}}

\newcommand{\be}{\begin{equation}}
\newcommand{\ee}{\end{equation}}
\newcommand{\bpm}{\begin{pmatrix}}
\newcommand{\epm}{\end{pmatrix}}

\numberwithin{equation}{section}

\newtheorem{theorem}{Theorem}
\newtheorem{lemma}{Lemma}
\newtheorem{proposition}{Proposition}
\newtheorem{corollary}{Corollary}

\theoremstyle{remark}
\newtheorem{definition}{Definition}

\newtheorem{remark}{Remark}
\newtheorem{example}{Example}

\begin{document}

\author[R. Bryant]{Robert L. Bryant}
\address{Duke University Mathematics Department\\
         P.O. Box 90320\\
         Durham, NC 27708-0320}
\email{\href{mailto:bryant@math.duke.edu}{bryant@math.duke.edu}}
\urladdr{\href{http://www.math.duke.edu/~bryant}%
         {http://www.math.duke.edu/\lower3pt\hbox{\symbol{'176}}bryant}}

\title[Gradient K\"ahler Ricci Solitons]
      {Gradient K\"ahler Ricci Solitons}

\date{August 02, 2004}

\begin{abstract}
Some observations about the local and
global generality of gradient K\"ahler Ricci solitons
are made, including the existence of a canonically
associated holomorphic volume form and vector
field, the local generality of solutions with
a prescribed holomorphic volume form and vector field,
and the existence of Poincar\'e coordinates in the 
case that the Ricci curvature is positive and the
vector field has a fixed point.
\end{abstract}



\thanks{
Thanks to Duke University for its support via
a research grant, to the NSF for its support
via DMS-8905207 and DMS-0103884 and to Columbia University 
for its support via an Eilenberg Visiting Professorship.  
\hfill\break
\hspace*{\parindent} 
This is Version~$2.0$ of RKSolitons.tex
}

\maketitle

\setcounter{tocdepth}{2}
\tableofcontents

\section{Introduction and Summary}\label{sec: intro}

This article concerns the local and global geometry 
of gradient K\"ahler Ricci solitons, i.e., K\"ahler
metrics~$g$ on a complex $n$-manifold~$M$ that admit 
a \emph{Ricci potential}, i.e., a function~$f$ 
such that~$\Ric(g) = \nabla^2 f$ (where~$\nabla$ 
denotes the Levi-Civita connection of~$M$.  

These metrics arise as limiting metrics in the
study of the Ricci flow~$g_t = -2\Ric(g)$ applied 
to K\"ahler metrics.  Under the Ricci flow, 
a gradient K\"ahler Ricci soliton~$g_0$ evolves 
by flowing under the vector field~$\nabla f$, i.e.,
\be
g(t) = {\exp_{(-t\nabla f)}}^*\bigl(g_0\bigr).
\ee
In particular, if the flow of~$\nabla f$ is complete,
then the Ricci flow with initial value~$g_0$ exists
for all time.

The reader who wants more background on these metrics 
might consult the references and survey articles
\cite{MR99g:53042,MR99j:53044,MR04f:32031}.  
The references \cite{math.DG/0310198,MR01g:53117,MR2001h:32040} 
contain further important work in the area and will
be cited further below.

\subsection{Basic facts}
Unless the metric~$g$ admits flat factors,
the equation~$\Ric(g) = \nabla^2 f$ determines~$f$
up to an additive constant and it does no harm to 
fix a choice of~$f$ for the discussion.  For simplicity,
it does no harm to assume that~$g$ has no
(local) flat factors and so this will frequently be 
done.  Also, the Ricci-flat case (aka the Calabi-Yau case), 
in which~$\Ric(g)=0$, is a special case that is usually
treated by different methods, so it will usually be
assumed that~$\Ric(g)\not=0$.  (Indeed, most of the 
latter part of this article will focus on the case
in which~$\Ric(g)>0$).

\subsubsection{The associated holomorphic vector field~$Z$}
One of the earliest observations~\cite{MR98a:53058} 
made about gradient K\"ahler Ricci solitons
is that the vector field~$\nabla f$
is the real part of a holomorphic vector field and
that, moreover,~$J(\nabla f)$ is a Killing field
for~$g$.  In this article, I will take~$Z =  
\frac12\bigl(\nabla f - \iC J(\nabla f)\bigr)$ to be
the holomorphic vector field \emph{associated} to~$g$.

\subsubsection{The holomorphic volume form~$\Upsilon$}
In the Ricci-flat case, at least when~$M$ is simply
connected, it is well-known that there is a $g$-parallel
holomorphic volume form~$\Upsilon$, i.e., one which satisfies
the condition that~$\iC^{n^2}2^{-n}\,\Upsilon\w\overline{\Upsilon}$
is the real volume form determined by~$g$ and the $J$-orientation.

In~\S\ref{ssec: assocholvolform}, I note that, for
any gradient K\"ahler Ricci soliton~$g$ with Ricci potential~$f$ 
defined on a simply connected~$M$, there is a holomorphic volume
form~$\Upsilon$ (unique up to a constant multiple of modulus~$1$)
such that $\iC^{n^2}2^{-n}\,\E^{-f}\,\Upsilon\w\overline{\Upsilon}$
is the real volume form determined by~$g$ and the $J$-orientation.
Of course,~$\Upsilon$ is not $g$-parallel (unless~$g$ is Ricci-flat)
but satisfies~$\nabla\Upsilon = \frac12\,\p f\otimes\Upsilon$.

This leads to a notion of \emph{special} coordinate charts 
for~$(g,f)$ i.e., coordinate charts~$(U,z)$ such that the 
associated coordinate volume form~$\d z = \d z^1\w\cdots\w\d z^n$ 
is the restriction of~$\Upsilon$ to~$U$.  In such coordinate
charts, several of the usual formulae simplify for 
gradient K\"ahler Ricci solitons.

\subsubsection{The $\Upsilon$-divergence of~$Z$}
Given a vector field and and volume form, the divergence
of the vector field with respect to the volume form is 
well defined.  It turns out to be useful to consider this
quantity for~$Z$ and~$\Upsilon$.  The divergence in this
case is the (necessarily holomorphic) function~$h$ that
satisfies~$\Lie_Z\Upsilon = h\,\Upsilon$.  

By general principles, the scalar function~$h$ must be 
expressible in terms of the first and second derivatives 
of~$f$.  Explicit computation 
(Proposition~\ref{prop: RnormZ2reln}) yields
\be
2h = \tr_g\bigl(\nabla^2 f) + |\nabla f|^2 = R(g) + |\nabla f|^2,
\ee
where~$R(g) = \tr_g\bigl(\Ric(g)\bigr)$ is the scalar curvature
of~$g$. In particular, $h$ is real-valued and therefore constant.
Now, the constancy of~$R(g) + |\nabla f|^2$ had already been
noted and utilized by Hamilton and Cao~\cite{MR01g:53117}.  
However, its interpretation as a holomorphic divergence 
seems to be new.

\subsection{Generality}
An interesting question is:  How many gradient K\"ahler
Ricci solitons are there?  Of course, this rather vague
question can be sharpened in several ways.  

The point of view adopted in this article is to start
with a complex $n$-manifold~$M$ already endowed with
a holomorphic volume form~$\Upsilon$ and a holomorphic
vector field~$Z$ and ask how many gradient K\"ahler solitons
on~$M$ there might be (locally or globally) that have~$Z$
and~$\Upsilon$ as their associated holomorphic data.

An obvious necessary condition is that the divergence~$h$
of~$Z$ with respect to~$\Upsilon$ must be a real constant.

\subsubsection{Nonsingular extension}
Away from the singularities (i.e., zeroes) of~$Z$, this
divergence condition turns out to be locally sufficient.

More precisely, I show (see Theorem~\ref{thm: locgen}) 
that if~$H\subset M$ is an embedded complex hypersurface 
that is transverse at each of its points to~$Z$, and~$g_0$
and~$f_0$ are, respectively, a real-analytic K\"ahler metric
and function on~$H$, then there is an open neighborhood~$U$
of~$H$ in~$M$ on which there exists a gradient K\"ahler
Ricci soliton~$g$ with potential~$f$ whose associated
holomorphic quantities are~$Z$ and~$\Upsilon$ and such 
that~$g$ and~$f$ pull back to~$H$ to become~$g_0$ and~$f_0$.
The pair~$(g,f)$ is essentially uniquely specified by
these conditions.  The real-analyticity of the `initial
data' $g_0$ and~$f_0$ is necessary in order for an
extention to exist since any gradient K\"ahler Ricci soliton 
is real-analytic anyway (see Remark~\ref{rem: realanalsolns}).

Roughly speaking, this result shows that, away from 
singular points of~$Z$, the local solitons~$g$ with associated
holomorphic data~$(Z,\Upsilon)$ depend on two arbitrary
(real-analytic) functions of $2n{-}2$ variables.

\subsubsection{Singular existence}
The existence of (local) gradient K\"ahler solitons
in a neighborhood of a singularity~$p$ of~$Z$ is both
more subtle and more interesting.

Even if the divergence of~$Z$ with respect to~$\Upsilon$
is a real constant, it is not true in general that a 
gradient K\"ahler Ricci solition with~$Z$ and~$\Upsilon$
as associated holomorphic data exists in a neighborhood 
of such a~$p$.

I show (Proposition~\ref{prop: Zlin})
that a necessary condition is that there exist
$p$-centered holomorphic coordinates~$z = (z^i)$ 
on a $p$-neighborhood~$U\subset M$ and real 
numbers~$h_1,\ldots,h_n$ such that, on~$U$,
\be
Z = h_1\,z^1\,\frac{\p\hfill}{\p z^1}
    + \cdots + h_n\,z^n\,\frac{\p\hfill}{\p z^n}.
\ee
In other words, $Z$ must be holomorphically linearizable,
with real eigenvalues.%
\footnote{Of course, it is by no means
true that every holomorphic vector field is 
holomorphically linearizable at each of its singular points.}

In such a case, if~$\Lie_Z\Upsilon= h\Upsilon$ 
where~$h$ is a constant, then~$h = h_1+\cdots+h_n$. 
I show (Proposition~\ref{prop: ZlinUpsconst}) that, 
moreover, in this case, one can always choose $Z$-linearizing
coordinates as above so that~$\Upsilon = \d z^1\w\cdots\w\d z^n$.

Thus, the possible local singular pairs~$(Z,\Upsilon)$ 
that can be associated to a gradient K\"ahler Ricci
soliton are, up to biholomorphism, parametrized by
$n$ real constants.

Using this normal form, one then observes that, 
by taking products of solitons of dimension~$1$,
any set of real constants~$(h_1,\ldots,h_n)$
can occur (see Remark~\ref{rem: singexistexamples}).
Since, for any gradient K\"ahler Ricci soliton~$g$
with associated holomorphic data~$(Z,\Upsilon)$, 
the formula~$\Ric(g) = \Lie_{\text{Re}(Z)} g$ holds, 
it follows that if~$g$ is such a K\"ahler
Ricci soliton defined on a neighborhood of a point~$p$
with~$Z(p)=0$, then~$h_1,\ldots,h_n$ are the eigenvalues
(each of even multiplicity) 
of~$\Ric(g)$ with respect to~$g$ at~$p$.

However, this does not fully answer the question of
how `general' the solitons are in a neighborhood of
such a~$p$.  In fact, this very subtly depends
on the numbers~$h_i$.  For example, if the $h_i\in\bbR$
are linearly independent over~$\bbQ$, then any 
gradient K\"ahler Ricci soliton~$g$ with associated
data~$(Z,\Upsilon)$ defined on a neighborhood of~$p$ 
must be invariant under the compact $n$-torus action
generated by the closure of the flow of the imaginary
part of~$Z$.  This puts severe restrictions on the
possibilities for such solitons.  

At the conclusion of~Section~\S\ref{sec: potentials}, 
I discuss the local generality problem near a 
singular point of~$Z$ and explain 
how it can best be viewed as an elliptic boundary 
value problem of a certain type, but do not go into
any further detail.  A fuller discussion of this
case may perhaps be undertaken at a later date.

\subsection{The positive case}
In Section~\S\ref{sec: poincare}, I turn to an
interesting special case:  The case where~$g$ is
complete, the Ricci curvature is positive, and
the scalar curvature~$R(g)$ attains its maximum
at some (necessarily unique) point~$p\in M$.

This case has been studied before by Cao and 
Hamilton~\cite{MR01g:53117}, who proved that
this point~$p$ is a minimum of the Ricci potential~$f$,
that~$f$ is a proper plurisubharmonic exhaustion function
on~$M$ (which is therefore Stein), and that, moreover,
the Killing field~$J(\nabla f)$ has a periodic orbit
on `many' of its level sets.

For simplicity, the Ricci potential~$f$ will be
be normalized so that~$f(p)=0$, so that~$f$ is
positive away from~$p$.

I show (Theorem~\ref{thm: poincarecoords}) that under these
assumptions there exist global~$Z$-linearizing coordinates
$z = (z^i):M\to\bbC^n$, so that~$M$ 
is biholomorphic to~$\bbC^n$ (which generalizes 
an earlier result of Chau and Tam~\cite{math.DG/0310198}).%
\footnote{On 27 July 2004, about 12 hours before the first version 
of this article was posted on the arXiv, Chau and Tam posted 
the first version their article~\texttt{math.DG/0407449} 
in which they prove, 
under the same hypotheses as in Theorem~\ref{thm: poincarecoords},
that~$M$ is biholomorphic to~$\bbC^n$.  I saw their posting
just before I posted this article.  Their method is different
and does not produce $Z$-linearizing coordinates, but has the
advantage that it applies in the case of expanding solitons.  
In the second (much shortened) version of their article, 
posted on 02 August, 2004, they deduce their biholomorphism
result from already-known results about automorphisms of
complex manifolds.}
Moreover, as a consequence, it follows that every positive 
level set of~$f$ has at least~$n$ periodic orbits, 
a considerable sharpening of Cao and Hamilton's original results.

This global coordinate system has several other applications.

For example, I show that there is a K\"ahler potential~$\phi$ 
for~$g$ that is invariant under the flow of~$J(\nabla f)$
and that this potential is unique up to an additive constant.
(Which can be normalized away by requiring that~$\phi(p)=0$.)

As another application, I show how to normalize the 
choice of~$Z$-linearizing holomorphic coordinates up 
to an ambiguity that lies in a compact subgroup of~$\Un(n)$.
This makes the function~$|z|$ well-defined on~$M$, so it
is available for estimates.  

As an illustration of such use,
I show that there are positive constants~$r$
and~$a_1$, $a_2$, $b_1$, $b_2$, $c_1$, and~$c_2$ 
such that, whenever~$|z|\ge r$,
\be
\begin{aligned}
a_1\,\log|z| \le f(&z) \le a_2\,\log|z|\,,\\
b_1\,\log|z| \le d(z&,0) \le b_2\,\log|z|\,,\\
c_1\,\bigl(\log|z|\bigr)^2 \le \phi(&z) \le c_2\,\bigl(\log|z|\bigr)^2\,.
\end{aligned}
\ee  
I also give some bounds for~$a_1$ and~$a_2$.
Perhaps these will be useful in further work.

\subsection{The toric case}

This section studies the geometry of 
the reduced equation in the case when 
a gradient K\"aher Ricci soliton~$g$
defined on a neighborhood of~$0\in\bbC^n$
has toric symmetry, i.e., is invariant under
the action of~$\bbT^n$, the diagonal subgroup
of~$\Un(n)$.  This may seem specialized,
but, for example, if the associated holomorphic
vector field is~$Z_\hs$ where~$\hs = (h_1,\ldots,h_n)$
and the real numbers $h_1,\ldots,h_n$ have the
`generic' property of being linearly 
independent over~$\bbQ$, then $g$ has toric symmetry.
Thus, metrics with toric symmetry are
the rule when~$Z$ has a `generic' singularity.

I first derive the equation satisfied by the
reduced potential, which turns out to be a singular
Monge-Amp\'ere equation.  (The singularities are,
of course, related to the singular orbits of the
$\bbT^n$-action.)  I then show that, nevertheless,
this singular equation has good regularity and
its singular initial value problem is well-posed
in the sense of G\`erard and Tahara~\cite{MR2001c:35056}.

As a consequence (Corollary~\ref{cor: singIVP}), 
it follows that, for any~$\hs\in\bbR^n$, 
any real-analytic $\bbT^{n-1}$-invariant 
K\"ahler metric on a neighborhood of~$0\in\bbC^{n-1}$ 
is the restriction to~$\bbC^{n-1}$ of an essentially
unique toric gradient K\"ahler Ricci soliton on
an open subset of~$\bbC^n$ with
associated holomorphic vector field~$Z = Z_\hs$ 
and associated holomorphic volume form~$\Upsilon = \d z$.
In particular, it follows that,
in a sense made precise in that section, 
the toric gradient K\"ahler Ricci solitons on~$\bbC^n$ 
depend on one `arbitrary' real-analytic function 
of~$(n{-}1)$ (real) variables.

Next, I show that the reduced (singular Monge-Amp\`ere)
equation is of Euler-Lagrange type, at least,
away from its singular locus, and discuss
some of its conservation laws via an application
of Noether's Theorem.  (This is in contrast
to the unreduced soliton equation, which is not 
variational).

\subsection{Acknowledgement}\label{ssec: thanks}
This work is mostly based on notes written 
after a conversation with Richard Hamilton 
during a visit he made to Duke University 
in November~1991.  Section~\ref{sec: poincare} 
is recent, having been written in April and May of~2004
after further conversations with Hamilton during
a semester I spent at Columbia University.

It is a pleasure thank Hamilton for his interest
and to thank Columbia University for its hospitality.

\section{Associated Holomorphic Quantities}\label{sec: holoquant}

In this section, constructions of some holomorphic
quantities associated to a gradient K\"ahler Ricci soliton~$g$
on a complex $n$-manifold~$M^n$ with Ricci potential~$f$
will be described.

\subsection{Preliminaries}
In order to avoid confusion because of various different
conventions in the literature, I will collect the notations,
conventions, and normalizations to be used in this article.

\subsubsection{Tensors and inner products}
Factors of~$2$ are sometimes troubling and confusing
in K\"ahler geometry.  

For~$a$ and~$b$ in a vector space~$V$, I will use
the conventions~$a\circ b = \frac12(a\otimes b + b \otimes a)$
and~$a\w b = a\otimes b - b\otimes a$.  In particular,
$a\otimes b = a\circ b + \frac12\,a\w b$.

A real-valued inner product~$\la,\ra$ on a real vector 
space~$V$ can be extended to~$V^\bbC = \bbC\otimes V$ 
in several different ways.  A natural way is to 
extend it as an Hermitian form, i.e., so that
\be
\la v_1 + \iC v_2, w_1 + \iC w_2\ra
= \bigl(\la v_1, w_1\ra + \la v_2, w_2\ra\bigr)
 +\iC \bigl(\la v_2, w_1\ra - \la v_1, w_2\ra\bigr)
\ee
and that is the convention to be adopted here.

If the real vector space~$V$ has a complex 
structure~$J:V\to V$, then $V^\bbC = V^{1,0}\oplus V^{0,1}$
where~$V^{1,0}$ is the $+\iC$-eigenspace of~$J$ extended
complex linearly to~$V^\bbC$ while~$V^{0,1}$ is the 
$(-\iC)$-eigenspace of~$J$.  It is common practice
to identify~$v\in V$ with~$v^{1,0} = v - \iC\,Jv\in V^{1,0}$,
but some care must be taken with this.

For example, an inner product~$\la,\ra$ on~$V$
is \emph{compatible} with~$J$ if~$\la Jv, Jw\ra = \la v, w\ra$
for all~$v,w\in V$.  Note the identity
\be
\la v^{1,0}, v^{1,0}\ra = 2\la v,v\ra.
\ee

For any $J$-compatible inner product~$\la,\ra$ on~$V$ 
(or equivalently, quadratic form)
there is an \emph{associated} $2$-form~$\eta$ defined by
\be
\eta(v,w) = \la Jv,w\ra\,.
\ee

\subsubsection{Coordinate expressions and the Ricci form}
Let $z=(z^i):U\to\bbC^n$ be a holomorphic coordinate chart on 
an open set $U\subset M$.  
The metric~$g$ restricted to~$U$ can be expressed in
the form
\be
g =  g_{i\bj}\,\d z^i{\circ}\d\zbar^j
\ee
for some functions~$g_{i\bj} = \overline{g_{j\bi}}$ on~$U$.
The associated K\"ahler form~$\Omega$ then has the 
coordinate expression
\be
\Omega = {\ts{\frac{\iC}{2}}}\,g_{i\bj}\>\d z^i\!\w \d\zbar^j\,.
\ee
Note that~$g_{i\bj}\,\d z^i{\otimes}\d\zbar^j = g - 2\iC\,\Omega$.

The Ricci tensor~${\rm Ric}(g)$ is $J$-compatible 
since~$g$ is K\"ahler, and hence has a coordinate 
expression~$\Ric(g) = R_{j\bk}\,\d z^j{\circ}\d\zbar^k$
where~$R_{j\bk} = \overline{R_{k\bj}}$. 
Its associated $2$-form~$\rho$ is computed by the formula
\be\label{eq: rhodef}
\rho = {\ts{\frac{\iC}{2}}}\,R_{i\bj}\>\d z^i\!\w \d\zbar^j
= -\iC\,\p{\overline\p}\,G
\ee
where
\be\label{eq: Gdef}
G=\log\,\det\bigl(g_{i\bj}\bigr).
\ee
While~$\rho$ is independent of the coordinate
chart used to compute it, the function~$G$ does depend
on the coordinate chart.

The scalar curvature~$R(g) = \tr_g\bigl(\Ric(g)\bigr)$
has the coordinate expression
\be
R(g) = 2 g^{i\bj} R_{i\bj}
\ee
and satisfies $R(g)\,\Omega^n = 2n\,{\ds\rho\w\Omega^{n-1}}$.

\subsubsection{The gradient K\"ahler Ricci soliton condition}
The following equivalent formulation of the gradient
K\"ahler Ricci soliton condition is well-known:

\begin{proposition}\label{prop: gKRSconditions}
A real-valued function~$f$ on~$M$ satisfies~$\Ric(g)=D^2f$ 
if and only if~$\rho=\iC\,\p{\overline\p}\,f$
and~$D^{0,2}f=0$.  This latter condition is
equivalent to the condition that the $g$-gradient of~$f$
be the real part of a holomorphic vector field on~$M$. \qed
\end{proposition}

\subsection{The associated holomorphic volume form}
\label{ssec: assocholvolform}
In this subsection, given a gradient K\"ahler
Ricci soliton~$g$ with Ricci potential~$f$ 
on a simply-connected complex $n$-manifold~$M$,
a holomorphic volume form on~$M$ (unique up to
a complex multiple of modulus~$1$) will be constructed.

\subsubsection{Existence of special coordinates}
The following result shows that there are coordinate
systems in which the Ricci potential is more closely
tied to the local coordinate quantities. 

\begin{proposition}\label{prop: speccoordexist}
If~$g$ is a gradient K\"ahler Ricci soliton on~$M$ with
Ricci potential~$f$, then~$M$ has an atlas of holomorphic 
charts~$(U,z)$ satisfying $\log\,\det\bigl(g_{i\bj}\bigr)=-f$.
\end{proposition}

\begin{proof} To begin, let~$(U,z)$ be any local
holomorphic coordinate chart on~$M$, with quantities~$g_{i\bj}$
and~$G$ defined as above.

Since~$f$ is a Ricci potential for~$g$, 
i.e., ~$\Ric(g)=D^2f$, it follows from~\eqref{eq: rhodef} 
and Proposition~\ref{prop: gKRSconditions} that
\be 
-\iC\,\p{\overline\p}\,G 
= \iC\,\p{\overline\p}\,f\,.
\ee
Thus,~$f+G$ is pluriharmonic. 
Assuming further that the domain~$U$ of the coordinate system~$z$
is simply connected, there exists a holomorphic function~$p$
on~$U$ so that 
\be
f = -G + p + {\bar p}.
\ee
Now let $w$ be any other local coordinate system 
on the same simply connected domain~$U$ in~$M$ and write 
\be
\Omega = {\ts\frac{\iC}2}\,h_{i\bj}\,\d w^i\!\w \d\wbar^j.
\ee
Then $H = \log\,\det\bigl(h_{i\bj}\bigr)$ is of the 
form 
\be
H = G + J + {\overline J}
\ee
where~$J$ is the log-determinant of the Jacobian matrix of the 
change of variables from $z$ to $w$, i.e.,
\be
\d z^1\w \d z^2\w\cdots\w \d z^n 
= e^J\,\d w^1\w\d w^2\w\cdots\w\d w^n\,.
\ee

It follows that every point of~$U$ has an open neighborhood~$V$ 
on which there exists a coordinate chart~$w$ for which $-H=f$, 
the Ricci potential.
\end{proof}

\begin{definition}[Special coordinates]
Let~$g$ be a gradient K\"ahler Ricci soliton on~$M$
with Ricci potential~$f$.
A coordinate chart $(U,z)$ for which  
$\log\det(g_{i\bj})=-f$ will be said to be \emph{special} 
for~$(g,f)$.
\end{definition}

\begin{remark}[The volume form in special coordinates]
A coordinate chart~$(U,z)$ is special for~$(g,f)$
if and only if the volume form of~$g$ satisfies
\be
\d\text{vol}_g 
= \frac1{n!}\,\Omega^n 
= \left(\frac{\iC^n}{2}\right)^n\,\E^{-f}\,\d z\w \d\zbar.
\ee
\end{remark}

\begin{theorem}[Existence of holomorphic volume forms]
\label{thm: holvolext}
Let~$M$ be a simply connected complex $n$-manifold
endowed with a gradient K\"ahler Ricci soliton~$g$
with associated K\"ahler form~$\Omega$ and a choice 
of Ricci potential~$f$.  Then there exists a holomorphic 
volume form~$\Upsilon$ on~$M$, unique up to muliplcation
by a complex number of modulus~$1$, with the property that
\be\label{eq: Holvoldef}
\d\text{vol}_g 
= \frac1{n!}\,\Omega^n 
= \left(\frac{\iC^n}{2}\right)^n\,
\E^{-f}\,\Upsilon\w \overline{\Upsilon}.
\ee
\end{theorem}

\begin{proof}
For any two $(g,f)$-special coordinate charts $z$ and $w$ on the
same domain~$U$, the ratio of their corresponding holomorphic 
volume forms is a constant of modulus~$1$.

The volume forms of special coordinate systems are thus the
sections of a flat connection~$\nabla_0$ on the canonical bundle of~$M$,
i.e., the bundle whose sections are the holomorphic volume forms on~$M$.  
Since~$M$ is simply connected, the canonical bundle of~$M$ 
has a global $\nabla_0$-flat section~$\Upsilon$ that is unique 
up to a multiplicative constant.

By construction, $\Upsilon$ satisfies~\eqref{eq: Holvoldef}.
Its uniqueness up to multiplication by a constant of modulus~$1$
is now evident.
\end{proof}

\begin{definition}[Associated holomorphic volume forms]
Given a gradient K\"ahler Ricci soliton~$g$ with Ricci potential~$f$,
a holomorphic volume form~$\Upsilon$ satisfying~\eqref{eq: Holvoldef} 
will be said to be \emph{associated} to the pair~$(g,f)$.
\end{definition}

\begin{remark}[Scaling effects on~$\Upsilon$]
Scaling a gradient K\"ahler Ricci soliton~$g$ by a constant
produces another gradient K\"ahler Ricci soliton and
adding a constant to~$f$ will produce another Ricci potential
for~$g$.

If~$\Upsilon$ is associated to~$(g,f)$, then, 
for any real constants~$\lambda>0$ and~$c$,
the $n$-form~$\lambda^n\E^c\,\Upsilon$ 
is associated to~$(\lambda^2\, g,\,f{+}2c)$.
\end{remark}

\subsection{The holomorphic flow}
Write the $g$-gradient of~$f$ as~$Z + {\bar Z}$ 
where~$Z$ is of type~$(1,0)$.
Thus, $Z = \frac12\bigl(\nabla f - \iC\,J(\nabla f)\bigr)$.

\subsubsection{The infinitesimal symmetry}
By the standard K\"ahler identities, 
$Z$ is the unique vector field of type~$(1,0)$ satisfying
\be\label{eq: fandZ}
{\bar\p}f = -\iC\,Z\lhk\Omega\,.
\ee
Writing~$Z = X - i\,Y = X - i\,JX$, it follows that, 
in addition to~$X$ being the one-half the gradient 
of~$f$, the vector field~$Y = JX$ is $\Omega$-Hamiltonian.
Thus, the flow of~$Y$ preserves~$\Omega$. 
 
Since~$Z$ is holomorphic by Proposition~\ref{prop: gKRSconditions}, 
the flow of~$Y$ also preserves the complex structure on~$M$.

Hence,~$Y$ must be a Killing vector field for the metric~$g$. 

Thus, a gradient K\"ahler Ricci soliton that is not Ricci-flat
always has a nontrivial infinitesimal symmetry.

\begin{proposition}
The singular locus of~$Z$ is a disjoint union of 
nonsingular complex submanifolds of~$M$, each of
which is totally geodesic in the metric~$g$. 
\end{proposition}

\begin{proof}
Since~$Z$ is holomorphic, its singular locus (i.e., the
locus where it vanishes) is a complex subvariety of~$M$.
However, since this locus is also the zero locus of~$Y 
= -\textrm{Im}(Z)$, which is a Killing field for~$g$,
this locus is a submanifold that is totally geodesic 
with respect to~$g$.  In particular, it must be smooth 
and hence nonsingular as a complex subvariety.
\end{proof}

\subsubsection{$Z$ in special coordinates}
Assume~$(U,z)$ is a special local coordinate system.  Since
\be
{\bar\p}G 
= g^{i\bj}\,{\frac{\p g_{i\bj}}{\p \zbar^k}}\,\d\zbar^k
= - {\bar\p} f \,,
\ee
the formula for~$Z$ in special coordinates is
\be\label{eq: Zformula}
Z =  Z^\ell {\frac{\p\hfil}{\p z^\ell}} = 
-\left(2g^{\ell\bk}g^{i\bj}{\frac{\p      
                 g_{i\bj}}{\p\zbar^k}}\right)
              {\frac{\p\hfil}{\p z^\ell}}\,.
\ee
Thus, the equations for a gradient K\"ahler Ricci 
soliton in special coordinates are that the functions $Z^\ell$  
defined by~\eqref{eq: Zformula} be holomorphic.

In fact, the expression in~\eqref{eq: Zformula} can be simplified, 
since the closure of~$\Omega$ is equivalent to the equations
\be\label{eq: Omegaclosed}
 {\frac{\p g_{i\bj}}{\p \zbar^k}}
={\frac{\p g_{i\bk}}{\p \zbar^j}}.
\ee
Thus,
\be\label{eq: Zasgup}
Z^\ell  = 
  -2\,g^{\ell\bk}g^{i\bj}\,{\frac{\p g_{i\bj}}{\p\zbar^k}}
= -2\,g^{i\bj}g^{\ell\bk}\,{\frac{\p g_{i\bk}}{\p\zbar^j}}
=2\,g^{i\bj}g_{i\bk}\,{\frac{\p g^{\ell\bk}}{\p\zbar^j}}
= 2 \,{\frac{\p g^{\ell\bj}}{\p\zbar^j}},
\ee
where I have used the identity~$g^{i\bj}g_{i\bk}=\delta^\bj_\bk$
and the identity~$g_{i\bk} g^{\ell\bk}=\delta^\ell_i$ and its derivatives.

\subsubsection{The $\Upsilon$-divergence of~$Z$}
Since $Z$ is holomorphic, the Lie derivative of~$\Upsilon$ 
with respect to $Z$ must be of the form $h\,\Upsilon$ 
where $h$ is a holomorphic function on $M$ (usually called
the divergence of~$Z$ with respect to~$\Upsilon$).

Replacing $\Upsilon$ by $\lambda\Upsilon$ for any 
$\lambda\in\bbC^*$ will not affect the definition of~$h$, 
so the function~$h$ is intrinsic to the geometry of the soliton.
On general principle, it must be computable in terms of the
first and second covariant derivatives of~$f$, which leads to
the following interpretation of a result of Cao and Hamilton:

\begin{proposition}
\label{prop: RnormZ2reln}
The holomorphic function~$h$ is real-valued 
{\upshape(}and therefore constant{\upshape)}.  
Moreover, 
\be\label{eq: hRZ2 reln}
2h= R(g) + 2|Z|^2
\ee
where~$R(g)$ is the scalar curvature of~$g$ 
and~$|Z|^2$ is the squared $g$-norm of~$Z$.
\end{proposition}

\begin{proof}
In special coordinates, where~$\Upsilon=\d z^1\w\cdots\w\d z^n$, 
the function~$h$ has the expression 
\be
h  = {\frac{\p Z^\ell}{\p z^\ell}}.
\ee
Thus, by~\eqref{eq: Zasgup},
\be
h = 2 \,{\frac{\p g^{\ell\bj}}{\p z^\ell\,\p\zbar^j}},
\ee
which shows that the holomorphic function~$h$ is real-valued 
and therefore constant.  Moreover, since~$\rho = \iC\,\p{\bar\p}f$, it follows that
\be\label{eq: Ricciformula}
\begin{aligned}
\left(\ts\frac\iC2\right)\,R_{j\bk}\,\d z^j\w \d \zbar^k
= \rho &= \iC\,\p{\bar\p}f = \p(Z\lhk \Omega)\\
&= \left(\ts\frac\iC2\right)\,
        \p\left(g_{\ell\bk}Z^\ell\,\d \zbar^{k}\right)\\
&= \left(\ts\frac\iC2\right)\,
       \left(
        g_{\ell\bk}\frac{\p Z^\ell}{\p z^j}
        +Z^\ell\frac{\p g_{\ell\bk}}{\p z^j}\right)\,
              \d z^j\w \d \zbar^{k}.
\end{aligned}
\ee
In particular, in view of~\eqref{eq: Omegaclosed}
and~\eqref{eq: Zformula}, 
\be
\begin{aligned}
R(g) &= 2 g^{j\bk}R_{j\bk}
= 2 g^{j\bk}\left(g_{\ell\bk}\frac{\p Z^\ell}{\p z^j}
        +Z^\ell\frac{\p g_{\ell\bk}}{\p z^j}\right)
= 2h + 2g^{j\bk}Z^\ell\frac{\p g_{\ell\bk}}{\p z^j}\\
&= 2h + 2Z^\ell\,g^{j\bk}\frac{\p g_{j\bk}}{\p z^\ell}
= 2h - g_{\ell\bk} Z^\ell {\bar Z}^k = 2h - 2|Z|^2,
\end{aligned}
\ee
as claimed.
\end{proof}

\begin{remark}[Interpretations]
\label{rem: CaoHamconslaw}
It was Cao and Hamilton~\cite[Lemma~4.1]{MR01g:53117} 
who first observed that the quantity~$R(g)+|\nabla f|^2$
is constant for a (steady) gradient K\"ahler Ricci soliton.
Since $Z = \frac12\bigl(\nabla f - \iC\,J(\nabla f)\bigr)$, 
one has~$2|Z|^2 = |\nabla f|^2$, so their expression 
is the right hand side of~\eqref{eq: hRZ2 reln}. 

The interpretation of $R(g)+|\nabla f|^2$ as the
$\Upsilon$-divergence of~$Z$ seems to be new.  

In a sense, this constancy can be regarded as a sort
of conservation law for the Ricci flow.  Note that,
since~$\Delta f = R(g)$, this relation
is equivalent to the equation~$\Delta_g (\E^f) = 2h\,\E^f$.
\end{remark}

\subsection{Examples}
The associated holomorphic objects constructed
so far make it possible to simplify somewhat the usual
treatment of the known explicit examples. The following
examples will be useful in later discussions in this
article.

\begin{example}[The one-dimensional case: Hamilton's cigar]
\label{ex: Hamsoliton}
Suppose that~$M$ is a Riemann surface.  Then~$\Upsilon$ is
a nowhere vanishing 1-form on~$M$ and~$Z$ is a holomorphic
vector field on~$M$ that 
satisfies~$\d\bigl(\Upsilon(Z)\bigr) = h\,\Upsilon$,
where~$h$ is a constant.
There are essentially two cases to consider.  

First, suppose that $h=0$.  
Then~$\Upsilon(Z)$ is a constant, say~$\Upsilon(Z) = c$.

If~$c=0$, then~$Z$ is identically zero, and, 
from~\eqref{eq: Zasgup} it follows that, in
special coordinates~$z = (z^1)$ the real-valued
function~$g^{1\bar1}$ is constant.  In particular, in
special coordinates~$g = g_{1\bar1}|\d z^1|^2$, so~$g$
is flat.

If~$c\not=0$, then~$Z$ is nowhere vanishing and, after
adjusting~$\Upsilon$ and the special coordinate system 
by a constant multiple, it can be assumed that~$c=2$, 
i.e., that $\Upsilon = \d z^1$ and $Z = 2\,\p/\p z^1$.  
Then~\eqref{eq: Zasgup} implies that~$g^{1\bar1}
= z^1 + \zbar^1 + C$ for some constant~$C$. By adding
a constant to~$z^1$, it can be assumed that~$C=0$, so
it follows that, in this coordinate system
\be\label{eq: gformh=0cnot0}
g = \frac{|\d z^1|^2}{(z^1 + \zbar^1)}.
\ee
Since~$M$ is supposed to be simply connected, one can 
take~$z^1$ to be globally defined.
Thus~$M$ is immersed into the right half-plane in~$\bbC$
in such a way that~$g$ is the pullback of the conformal 
metric defined by~\eqref{eq: gformh=0cnot0}.  Of course,
this metric is not complete, even on the entire right
half-plane.

	Second, assume that~$h$ is not zero.  Then
$\Upsilon(Z)$ is a holomorphic function on~$M$ that
has nowhere vanishing differential.  Write~$\Upsilon(Z)=h z^1$
for some (globally defined) holomorphic immersion~$z^1:M\to\bbC$.
Then, by construction,~$\Upsilon = \d z^1$ 
and~$Z = hz^1\,\p/\p z^1$.  
By~\eqref{eq: Zasgup}, it follows that
\be
g^{1\bar1} = {\ts\frac12}(c + h\,|z^1|^2)
\ee
for some constant~$c$, so $z^1(M)\subset\bbC$ must lie
in the open set~$U$ in the~$w$-plane on which $c+h|w|^2 >0$.
In fact, $g$ must be the pullback under~$z^1:M\to U\subset\bbC$ of
the metric
\be\label{eq: Hamiltoncigar}
\frac{2\,|\d w|^2}{c+h\,|w|^2}.
\ee
This metric on the domain~$U\subset\bbC$
is not complete unless both~$c$ and~$h$ are nonnegative
and it is flat unless both~$c$ and~$h$ are positive.  In
this latter case, this metric is simply Hamilton's 
`cigar' soliton~\cite{MR89i:53029}.

Consequently, in dimension~$1$, the only complete gradient
K\"ahler Ricci solitons are either flat or one of Hamilton's 
`cigar' solitons (which are all homothetic to a single 
example).

Note that, under the Ricci flow~$g_t = -2\Ric(g)$, the
metric~\eqref{eq: Hamiltoncigar} evolves as
\be
g(t) = \frac{2\,|\d w|^2}{\E^{2ht}c+h\,|w|^2} 
= \frac{2\,\bigl|\d (\E^{-ht}w)\bigr|^2}{c+h\,\bigl|\E^{-ht}w\bigr|^2} 
= \Phi(-t)^*(g_0)
\ee
where~$\Phi(t)(w) = e^{ht}w$ is the flow of twice the
real part of~$Z = hw\,\p/\p w$. 
\end{example}

\begin{example}[Products]\label{ex: products}
By taking products of the $1$-dimensional examples, one
can construct a family of complete examples on~$\bbC^n$:
Let~$h_1,\ldots,h_n$ and $c_1,\ldots,c_n$ be positive
real numbers and consider the metric on~$\bbC^n$ defined
by
\be\label{eq: productgKRs}
g = \sum_{k=1}^n \frac{2\,|\d w^k|^2}{\left(c_k+h_k\,|w^k|^2\right)}.
\ee
This is, of course, a gradient K\"ahler Ricci soliton, 
with associated holomorphic volume form and vector field
\be
\Upsilon = \d w^1\w\d w^2\w\cdots\w\d w^n,
\qquad
Z = \sum_{k=1}^n h_k\,w^k\,\frac{\p\hfill}{\p w^k}.
\ee  
The Ricci curvature is
\be
\Ric(g) = 
\sum_{k=1}^n \frac{2c_kh_k\,|\d w^k|^2}{\left(c_k+h_k\,|w^k|^2\right)^2}.
\ee 

Although these product examples are trivial generalizations
of Hamilton's cigar soliton, they will be useful in observations
to be made below. 

Also, note that, even if the~$h_k$ are not positive,
as long as the $c_k$ are positive, the 
formula~\eqref{eq: productgKRs} defines a not-necessarily-complete
gradient K\"ahler Ricci soliton on the polycylinder defined
by the inequalities~$c_k+h_k|w^k|^2>0$.
\end{example}

\begin{example}[Cao's Soliton]\label{ex: Caosoliton}
One more case of an easily constructed example is the
gradient K\"ahler Ricci soliton metric~$g$ on~$\bbC^n$ 
that is invariant under~$\Un(n)$, 
discovered by H.-D. Cao~\cite{MR98a:53058}.  The
form of this metric can be derived as follows:  

Suppose that such a metric~$g$ is given on~$\bbC^n$.  
(One could do this analysis on any~$\Un(n)$-invariant 
domain in~$\bbC^n$, and Cao does this, but I will not pursue 
this more general case further here.)
The group~$\Un(n)$ must preserve the associated
holomorphic volume form~$\Upsilon$ up to a constant multiple 
and this implies that~$\Upsilon$ must be a constant multiple
of the standard volume form~$\d z^1\w\cdots\w\d z^n$. 
Since~$\Upsilon$ is only determined up to a constant multiple
anyway, there is no loss of generality in assuming that
$\Upsilon = \d z^1\w\cdots\w\d z^n$.  Furthermore, the 
vector field~$Z$ must also be invariant under~$\Un(n)$,
which implies that~$Z$ must be a multiple of the radial vector
field.  Since~$\d\bigl(Z\lhk\Upsilon) = h\,\Upsilon$ where~$h$
is real, it follows that 
\be
Z = h\,\sum_{k=1}^n \,z^k\,\frac{\p\hfill}{\p z^k}\,.
\ee

Now, the condition that~$g$ be rotationally invariant with
associated K\"ahler form closed implies that
\be\label{eq: grotsymmasa}
g_{i\bj} = a(r)\delta_{ij} + a'(r)\,\bar z^i\,z^j
\ee
for some function~$a$ of~$r = |z^1|^2{+}\cdots{+}|z^n|^2$
that satisfies~$ra'(r)+a(r)>0$ and~$a(r)>0$ (when~$n>1$).
Thus~$G = \log\bigl(a(r)^{n-1}(ra'(r){+}a(r))\bigr)$ in 
this coordinate system.  Now, the identity~$G = -f$,
the equation~\eqref{eq: fandZ}, and the above formula for 
the coefficients of~$\Omega$ combine to yield 
\be
\bar\p G = \iC\,Z\lhk\Omega 
= -\frac{h}{2}\bigl(ra'(r){+}a(r)\bigr) \bar\p r
= -\frac{h}{2}\bar\p\bigl(ra(r)\bigr).
\ee
Supposing that~$n>1$ (since the $n=1$ case has already been
treated), it follows that~$G+\frac{h}{2}ra(r)$ must be constant,
i.e., that
\be\label{eq: asingODE}
a(r)^{n-1}\bigl(ra(r)\bigr)' e^{(h/2)ra(r)} = a(0)^n.
\ee
Upon scaling~$\Upsilon$ by a constant, it can be assumed
that~$a(0)=1$, so assume this from now on.  Also, one can
assume that~$h$ is nonzero since, otherwise, the solution
that is smooth at~$r=0$ is simply~$a(r)\equiv a(0) = 1$,
which gives the flat metric.

The \textsc{ode}~\eqref{eq: asingODE} for~$a$ is
singular at~$r=0$, so the existence of a smooth 
solution near~$r=0$ is not immediately apparent. 

Fortunately, \eqref{eq: asingODE} can be integrated by
quadrature:  Set~$b(r) = (h/2)ra(r)$ and note that 
\eqref{eq: asingODE} can be written in terms of~$b$ as
\be\label{eq: bsingODE}
b(r)^{n-1}e^{b(r)}\,b'(r) = (h/2)^nr^{n-1}.
\ee
Integrating both sides from~$0$ to~$r>0$ yields an equation
of the form
\be\label{eq: Fbreln}
(-1)^n(n{-}1)!\,e^{b(r)}
\left(e^{-b(r)}-\sum_{k=0}^{n-1}\frac{\bigl(-b(r)\bigr)^k}{k!}
\right)
= \left(\frac{h}2\right)^n\,\frac{r^n}n.
\ee
Set
\be
F(b) = (-1)^n(n{-}1)!\,e^b
\left(e^{-b}-\sum_{k=0}^{n-1}\frac{(-b)^k}{k!}
\right)
\simeq e^b\left(\frac{b^n}{n} - \frac{b^{n+1}}{n(n{+}1)} + \cdots\right).
\ee

Now,~$F$ has a power series of the form~$F(b) 
= \frac1n\,b^n(1+ \frac{n}{n+1}b + \cdots)$,
so~$F$ can be written in the form~$F(b) = \frac1{n} f(b)^n$
for an analytic function of the form~$f(b) = b(1+\frac1{n+1}b+\cdots)$.
The analytic function~$f$ is easily seen to satisfy~$f'(b)>0$
for all~$b$ and to satisfy the limits
\be
\lim_{b\to+\infty} f(b) = \infty
\qquad\text{and}\qquad
\lim_{b\to-\infty} f(b) = -{\root n \of {n!}}\,.
\ee 
In particular, $f$ maps~$\bbR$ diffeomorphically
onto~$\bigl(-{\root n \of {n!}},\infty\bigr)$ 
and is smoothly invertible.  Of course,~$f(0)=0$.

Since~\eqref{eq: Fbreln} is equivalent to~$f\bigl(b(r)\bigr)^n
= \bigl(\frac{h}2\,r\bigr)^n$, when~$h>0$ it can be solved
for~$r\ge0$ by setting~$b(r) = f^{-1}\bigl(\frac{h}2\,r\bigr)$,
yielding a unique real-analytic solution with a power series
of the form
\be
b(r) = \frac{h}2\,r - \frac{h^2}{4(n{+}1)}\,r^2 + \cdots\ {} .
\ee

Consequently, when~$h>0$, the solution~$b$ is defined for all~$r\ge0$
and is positive and strictly increasing on the half-line~$r\ge0$.
In particular, the function
\be
a(r) = \frac2{h}\,\frac{b(r)}{r}
= 1 - \frac{h}{2(n{+}1)}\,r + \cdots\ {} .
\ee
is a positive real-analytic solution of~\eqref{eq: asingODE} 
that is defined on the range~$0\le r < \infty$ and
satisfies~$ra'(r)+a(r) = b'(r) > 0$ on this range, so 
that the expression~\eqref{eq: grotsymmasa} defines a
gradient K\"ahler Ricci soliton on~$\bbC^n$.

An \textsc{ode} analysis of this solution (which 
Cao~\cite{MR98a:53058} does) shows that when~$h>0$ 
the resulting metric is complete on~$\bbC^n$ 
and has positive sectional curvature.

When~$h<0$, the solution~$b(r)$ 
only exists for~$r < -\frac2h\root n \of{n!}\,$. 
It is not difficult to see that the corresponding 
gradient K\"ahler Ricci soliton on a bounded ball 
in~$\bbC^n$ is inextendible and incomplete.
\end{example}

\section{Potentials and local generality}\label{sec: potentials}

In this section, the question of `how many' gradient K\"ahler Ricci
soliton metrics could give rise to specified holomorphic 
data~$(\Upsilon, Z)$ on a complex manifold~$M$ will be 
considered.  While this question is not easy to answer globally,
it is not so difficult to answer locally.

Thus, throughout this section, assume that a complex
$n$-manifold~$M$ is specified, together with a nonvanishing 
holomorphic volume form~$\Upsilon$ on~$M$ and a holomorphic
vector field~$Z$ on~$M$ such that~$\d\bigl(Z\lhk\Upsilon\bigr)
= h\,\Upsilon$ for some real constant~$h$.

\subsection{Local potentials}\label{ssec: locpot}
Suppose that~$U\subset M$ is an open subset on which
there exists a function~$\phi$ such that~$\Omega 
= {\frac\iC2}\,\p{\bar\p}\phi$ is a positive
definite $(1,1)$-form whose associated K\"ahler metric~$g$
is a gradient Ricci soliton with associated holomorphic
data~$\Upsilon$ and~$Z$ and Ricci potential~$f$.

By~\eqref{eq: fandZ}, 
\be
\begin{aligned}
2\bar\p f &= -2\iC Z\lhk\Omega 
= Z\lhk(\p\bar\p\phi) 
= -Z\lhk(\bar\p\p\phi) \\
&= -Z\lhk\bigl(\d(\p\phi)\bigr)
= -\Lie_Z(\p\phi) + \d\bigl(\p\phi(Z)\bigr)\\
&= \bar\p\bigl(\p\phi(Z)\bigr) 
   - \left(\Lie_Z(\p\phi) - \p\bigl(\Lie_Z(\phi)\bigr)\right)
\end{aligned}
\ee 
By decomposition into type, it follows that
\be\label{eq: 2f-pphiZhol}
\bar\p\bigl(2f - \p\phi(Z)\bigr)  = 0.
\ee
Consequently, $F = 2f - \p\phi(Z)= 2f - \d\phi(Z) $ 
is a holomorphic function on~$U$.

\subsection{Nonsingular extension problems}
\label{ssec: nonsingextprobs}
Suppose now that~$p\in U$ is not a singular point of~$Z$.
Then, by shrinking~$U$ if necessary, $F$ can be written 
in the form~$F = dH(Z)$ for some holomorphic function~$H$
on the $p$-neighborhood~$U$.  Replacing~$\phi$ by $\phi+ H + \bar H$,
gives a new potential for~$\Omega$ that satisfies the stronger
condition
\be\label{eq: phinormalized}
\p\phi(Z) = \d\phi(Z) = 2 f.
\ee
This function~$\phi$ is unique up to the addition of the
real part of a holomorphic function that is constant on
the orbits of~$Z$.

Of course, \eqref{eq: phinormalized} implies that~$\d\phi(Y)=0$,
i.e., that~$\phi$ is invariant under the flow of~$Y$, the
imaginary part of~$Z$.

\subsubsection{Local reduction to equations}
In local coordinates~$z = (z^i)$ for which~$\Upsilon
= \d z^1\w\cdots\w\d z^n$, one has~$f = -G$ 
so~$\phi$ satisfies the Monge-Amp\`ere equation%
\footnote{It is interesting to note that this
equation is not of Euler-Lagrange type, even locally, 
unless~$Z\equiv0$, i.e., the Ricci-flat case.
Of course, in the Ricci-flat case,
the variational nature of this equation is well-known.}
\be\label{eq: MAcond}
\det\left(\frac{\p^2\phi}{\p z^i\p\zbar^j}\right)
\E^{\frac12\d\phi(X)} = 1
\ee
as well as the equation
\be\label{eq: Yinvcond}
\d\phi(Y)=0.
\ee

Conversely, if~$\phi$ is a strictly pseudo-convex function defined 
on a $p$-neighborhood~$U$ that satisfies both~\eqref{eq: MAcond} 
and~\eqref{eq: Yinvcond}, then the K\"ahler metric~$g$ whose
associated K\"ahler form is~$\Omega 
= \frac\iC2\,\p\bar\p\phi$ is a gradient K\"ahler
Ricci soliton on~$U$ with associated holomorphic form~$\Upsilon$
and~holomorphic vector field~$Z$.

\begin{remark}[Real-analyticity of solitons]
\label{rem: realanalsolns}
Note that, because~\eqref{eq: MAcond} is a real-analytic elliptic 
equation for the strictly pseudo-convex function~$\phi$, 
it follows by elliptic regularity that~$\phi$ (and hence~$g$) 
is real-analytic as well.
\end{remark}

Now, \eqref{eq: MAcond} and~\eqref{eq: Yinvcond} are two
\textsc{pde} for~$\phi$, the first of second order and
the second of first order.  
While this is an overdetermined system, it is not difficult 
to show that it is involutive in Cartan's sense.  

In fact, an analysis along the lines of exterior differential 
systems leads to the following result as a proper formulation 
of a `Cauchy problem' for gradient K\"ahler Ricci solitons 
in the nonsingular case: 

\begin{theorem}[Nonsingular extensions]\label{thm: locgen} 
Let~$M^n$ be a complex $n$-manifold endowed with a holomorphic
volume form~$\Upsilon$ and a nonzero vector field~$Z$ satisfying
$\d(Z\lhk\Upsilon) = h\,\Upsilon$ for some real constant~$h$.

Let~$H^{n-1}\subset M$ be any embedded complex hypersurface 
that is transverse to~$Z$, let~$\Omega_0$ be any real-analytic K\"ahler
form on~$H$, and let~$f_0$ be any real-analytic function on~$H$.

Then there is an open $H$-neighborhood~$U\subset M$ on
which there exists a gradient K\"ahler Ricci soliton~$g$ with
associated K\"ahler form~$\Omega$, holomorphic volume form~$\Upsilon$,
holomorphic vector field~$Z$, and Ricci potential~$f$ that
satisfy%
\footnote{Notation:  If~$P\subset Q$ is a submanifold, 
and~$\psi$ is a differential form on~$Q$, I use~$P^*\phi$
to denote the pullback of~$\psi$ to~$P$.}
\be
H^*\Omega=\Omega_0,\qquad\text{and}\qquad H^*f = f_0.
\ee
Moreover,~$g$ is locally unique in the sense that any other
gradient K\"ahler Ricci soliton~$\tilde g$ defined on
an open~$H$-neighborhood~$\tilde U\subset M$ 
satisfying these initial conditions agrees with~$g$ 
on some open neighborhood of~$H$ in~$U\cap\tilde U$.
\end{theorem}  

\begin{proof}
The first step in the proof will be to construct a
special set of local `flow-box' coordinate charts adapted to the 
hypersurface~$H$, the holomorphic form~$\Upsilon$, and
the holomorphic vector field~$Z$.

To begin, note that, since, by hypothesis~$Z_p$ does not lie
in~$T_pH\subset T_pM$ for all~$p\in H$, 
the $(n{-}1)$-form~$Z\lhk\Upsilon$ is nonvanishing 
when pulled back to~$H$.

Let~$p\in H$ be fixed. Since~$H^*\bigl(Z\lhk\Upsilon\bigr)$
does not vanish at~$p$, there exist $p$-centered holomorphic 
coordinates~$w^2,\ldots,w^n$ on a $p$-neighborhood~$V$ in~$H$ 
such that~$V^*\bigl(Z\lhk\Upsilon\bigr) = \d w^2\w\cdots\w\d w^n$.

Since~$H$ is embedded in~$M$, there exists an open 
neighborhood~$U\subset M$ of~$V\subset H$ with the property 
that~$U\cap H = V$ and so that each complex integral 
curve~$C\subset M$ of~$Z$ that meets~$U$ does so in a 
connected open disk~$U\cap C$ that intersects~$H$ in 
a unique point.  

Consequently, there exist unique holomorphic 
functions~$z^2,\ldots,z^n$
on~$U$ satisfying $\d z^2(Z) = \cdots = \d z^n(Z) = 0$
and~$V^*(z^j) = w^j$ for~$2\le j\le n$.  Moreover, there
exists a unique function~$z^1$ on~$U$ with the property
that~$z^1$ vanishes on~$V = U\cap H$ and so that~$U^*\Upsilon
= \d z^1\w\d z^2\w\cdots\w\d z^n$.  Since the functions~$z^1,\ldots,z^n$
have independent differentials on~$U$, it follows that
by shrinking~$V$ (and hence~$U$) if necessary, it can be assumed 
that~$(U,z)$ is a $p$-centered holomorphic coordinate chart
whose image~$z(U)\subset\bbC^n$ is a polycylinder of the
form~$|z^i|<\rho^i$ for some~$\rho^1,\ldots,\rho^n>0$.
By shrinking~$\rho^1$ if necessary, it can be arranged that
$1{+}h\rho^1 >0$.

By construction,~$Z = F(z)\,\p/\p z^1$ for some holomorphic
function~$F$ defined on~$z(U)\subset\bbC^n$. 
Thus,~$U^*(Z\lhk\Upsilon) = F(z)\,\d z^2\w\cdots\w\d z^n$.
Since~$V^*(Z\lhk\Upsilon) = \d w^2\w\cdots\w\d w^n$, it follows
that~$F(0,w^2,\ldots,w^n) = 1$ for~$(0,w^2,\ldots,w^n)\in z(U)$.
Moreover, since~$\d\bigl(Z\lhk\Upsilon\bigr) = h\,\Upsilon$,
it follows that~$\p F/\p z^1 = h$.  Consequently, in these
coordinates~$Z = (1{+}h z^1)\,\p/\p z^1$.

Now write~$Z = X - \iC Y$, where~$X$ and~$Y$ are commuting 
real vector fields.  The integral curves of~$Y$ are transverse
to the hypersurface~$H$ and there exists a real 
hypersurface~$R\subset U$ that is the union of the integral curves
of~$Y$ in~$U$ that pass through~$V = U\cap H$.  The vector field~$X$
is everywhere transverse to~$R$ in~$U$.  

Now let~$\psi_0$ be a real-valued function on~$V$ such that
$V^*(\Omega_0) = \frac\iC2\,\p\bar\p\psi_0$.  Such an  
$\Omega_0$-potential~$\psi$ is unique up the the addition 
of the real part of a holomorphic function of~$w^2,\ldots, w^n$.
Extend~$\psi_0$ to a function~$\psi_1$ on~$R$ by making it constant
along the integral curves of~$Y$.  Similarly, extend~$V^*(f_0)$
to a function~$f_1$ on~$R$ by making it constant along the
integral curves of~$Y$.

Finally, consider the initial value problem for a 
function~$\phi$ on a neighborhood of~$R$ in~$U$ given by
the real-analytic \textsc{pde}
\be\label{eq: IVPe}
\det\left(\frac{\p^2\phi}{\p z^i\p\zbar^j}\right)
\E^{\frac12\d\phi(X)} = 1
\ee
subject to the real-analytic initial conditions
\be\label{eq: IVPc}
\begin{aligned}
\phi(z) &= \psi_1(z)\\
\Lie_X(\phi)(z) &= 2 f_1(z)
\end{aligned}
\qquad \text{for all~$z\in R\subset U$}.
\ee
It is easy to check that \eqref{eq: IVPe} and~\eqref{eq: IVPc}
constitutes a noncharacteristic Cauchy
problem.  Hence, by the Cauchy-Kovalewski Theorem, 
there exists an open neighborhood~$W\subset U$
containing~$R$ on which there exists a solution~$\phi$ 
to this problem.  

Now, the solution~$\phi$ produced by the
Cauchy-Kovalewski Theorem is real-analytic and strictly
pseudo-convex.  By uniqueness in the Cauchy-Kovalewski 
Theorem, $\phi$ is the unique real-analytic solution.  Since,
as has already been remarked, elliptic regularity implies that
any strictly pseudo-convex solution of~\eqref{eq: IVPe} must
be real-analytic, it follows that~$\phi$ is the unique 
solution of~\eqref{eq: IVPe} that satisfies~\eqref{eq: IVPc}.

By its very construction, the $(1,1)$-form
$\Omega = {\ts\frac\iC2}\p\bar\p\phi$ is then the K\"ahler form
of a gradient K\"ahler Ricci soliton metric on~$W\subset U$ that 
satisfies~$V^*\Omega = V^*\Omega_0$, that has~$W^*\Upsilon$
and~$W^*Z$ as the associated holomorphic volume form and
vector field, respectively, and has~$f = {\ts\frac12}\d\phi(X)$
as Ricci potential, which, of course, satisfies~$V^*f = V^*f_0$.

Now, if one replaces~$\psi$ by~$\psi + H + \bar H$ for
some holomorphic function~$H=H(w^2,\ldots,w^n)$ on~$V$, 
then one finds that the solution~$\phi$ 
is replaced by by~$\phi + H(z^2,\ldots,z^n) 
+ \overline{H(z^2,\ldots,z^n)}$, so that~$\Omega$ is unaffected.

The argument thus far has shown that every point~$p\in H$
has an open neighborhood~$U\subset M$ on which there exists
a gradient K\"ahler-Ricci soliton~$g_U$ with the desired 
extension properties and associated holomorphic data.  It has
also shown that this extension is locally unique.  Now a
standard patching argument shows that there exists an open
neighborhood~$U\subset M$ of the entire complex hypersurface~$H$
on which such an extension exists and is unique in the 
sense described in the statement of the theorem.
\end{proof}

\begin{remark}[Local generality]\label{rem: locgen}
Theorem~\ref{thm: locgen} essentially says that the local
gradient K\"ahler Ricci solitons depend on two real-analytic 
functions of~$2n{-}2$ variables, namely the potential 
functions~$\psi_0$ (which is assumed to be strictly
pseudo-convex but otherwise arbitrary) 
and~$f_0$ (which is arbitrary).  There is, of course, some
ambiguity in the choice of the holomorphic coordinates~$z^i$,
but this ambiguity turns out to depend on essentially~$n{-}2$ 
holomorphic functions of~$n{-}1$ holomorphic variables, 
which is negligible when compared with two arbitrary (real-analytic)
functions of $2n{-}2$ real variables.
\end{remark}

\subsection{Near singular points of~$Z$}\label{sssec: Zsingpts}
The situation near a singular point of~$Z$ is considerably
more delicate and interesting.

\subsubsection{Linear parts and linearizability}
Recall that, at a point~$p\in M$ where~$Z$ vanishes, there
is a well-defined linear map~$Z'_p:T_pM\to T_pM$ 
(often called `the linear part of~$Z$ at~$p$') defined
by setting~$Z'_p(v) = w$ if $w = [V,Z](p)$ for some (and
hence any) holomorphic vector field~$V$ defined near~$p$
and satisfying~$V(p)=v\in T_pM$.  

In local coordinates~$z = (z^i)$ centered on~$p$, if
\be
Z =  Z^j(z)\,\frac{\p\hfill}{\p z^j},
\ee
where, by assumption~$Z^j(0)=0$ for~$1\le j\le n$, then
\be
Z'_p\left(\,\frac{\p\hfill}{\p z^l}(p)\,\right)
= \frac{\p Z^j}{\p z^l}(0)\,\frac{\p\hfill}{\p z^j}(p).
\ee
The linear map~$Z'_p:T_pM\to T_pM$ has a Jordan
normal form and this is an important invariant of the
singularity.  In particular, the set of eigenvalues of~$Z'_p$
is well-defined.

\begin{proposition}\label{prop: Zprealeigens}
Let~$Z$ be the holomorphic vector field associated
to a gradient K\"ahler Ricci soliton~$g$ on~$M$.
At any singular point of~$Z$, the linear part~$Z'_p$
is diagonalizable, with all eigenvalues real.
\end{proposition}

\begin{proof}
If the data~$(\Upsilon, Z)$ is associated to a 
gradient K\"ahler Ricci soliton~$g$ in a neighborhood of
a singular point~$p$ of~$Z$,
then~\eqref{eq: Ricciformula} shows that, 
in special coordinates centered on~$p$, one has
\be\label{eq: Ricciformlinpart}
g^{i\bk}(0)\,R_{j\bk}(0) =  \frac{\p Z^i}{\p z^j}(0).
\ee
Because the matrices~$\bigl(g_{i\bj}(0)\bigr)$ 
and~$\bigl(R_{i\bj}(0)\bigr)$ are Hermitian symmetric
and~$\bigl(g_{i\bj}(0)\bigr)$ is positive definite, one
can choose the special coordinates so that~$\bigl(g_{i\bj}(0)\bigr)$
is a multiple of the identity matrix and~$\bigl(R_{i\bj}(0)\bigr)$
is diagonal.  
\end{proof} 

\begin{definition}\label{def: Zlin}
A holomorphic vector field~$Z$ on~$M$ is said to be
\emph{linearizable} near a singular point~$p$ if there
exist $p$-centered coordinates~$w = (w^i)$ on 
an open $p$-neighborhood~$W$ and constants~$a^i_j$ such
that, on~$W$, one has
\be
Z = a^i_j w^j\,\frac{\p\hfill}{\p w^i}\,.
\ee
The coordinates~$w = (w^i)$ are said to be \emph{linearizing}
or~\emph{Poincar\'e} coordinates for~$Z$ near~$p$.
\end{definition}

Not every holomorphic vector field is linearizable near
its singular points, even if the linear part at such a point
has all of its eigenvalues nonzero and distinct. 

\begin{example}[A nonlinearizable singular point]
\label{ex: notlinZ} The vector field 
\be
Z = z^1\,\frac{\p\hfill}{\p z^1} 
      + \bigl(2z^2 + (z^1)^2\bigr)\,\frac{\p\hfill}{\p z^2}
\ee
on~$\bbC^2$ is not linearizable at the origin, even though
its linear part there is diagonalizable with eigenvalues~$1$ and~$2$.

This nonlinearizability is perhaps most easily seen as follows:
The flow~$\Phi(t)$ of the vector field~$Z$ is 
\be
\Phi(t)(z^1,z^2) = \bigl(\E^t\,z^1,\,\E^{2t}(z^2 + (z^1)^2t)\bigr).
\ee
In particular~$\Phi(t+2\pi\iC)\not=\Phi(t)$, which would
be true if~$Z$ were holomorphically conjugate to the linear 
vector field
\be
Z'_{(0,0)} = z^1\,\frac{\p\hfill}{\p z^1} 
      + 2z^2\,\frac{\p\hfill}{\p z^2}.
\ee
\end{example}

This phenomenon, however, does not happen for singular points
of holomorphic vector fields associated to a gradient K\"ahler
Ricci soliton:

\begin{proposition}\label{prop: Zlin}
Let~$Z$ be a nonzero holomorphic vector field on
the complex $n$-manifold~$M$ that is associated 
to a gradient Kahler Ricci soliton~$g$.  Then $Z$ is 
linearizable at each of its singular points.
Moreover, the linear part of~$Z$ at a singular point
is diagonalizable and has all its eigenvalues real.
\end{proposition}

\begin{proof}
Let~$p\in M$ be a singular point of~$Z$.
The diagonalizability of the linear part of~$Z$ at
a singular point and the reality of the corresponding
eigenvalues has already been demonstrated, so all that
remains is to show that~$Z$ is linearizable near~$p$.

To do this, write~$Z = X - \iC Y$ where~$X$ and~$Y$ are,
as usual, real vector fields.  As has already been remarked,
the vector field~$Y$ is an infinitesimal isometry of~$g$.
In particular, the flow of~$Y$ is complete in the
geodesic ball~$B_r(p)$ for some~$r>0$ and is a $1$-parameter
group of isometries of the metric~$g$ restricted to~$B_r(p)$
that fixes the center~$p$.  
It follows that there is a compact, connected abelian 
subgroup~$\bbT\subset\Un(T_pM)$ whose Lie 
algebra is an abelian subalgebra~$\eut\subset\euu(T_pM)$
that contains~$Y'_p:T_pM\to T_pM$, the linearization of~$Y$ at~$p$
and is such that the $1$-parameter subgroup~$\exp(t Y'_p)$
is dense in~$\bbT$.

Let~$\Phi:\bbT\to\text{Isom}\bigl(B_r(p),g\bigr)$ be the
homomorphism induced by the exponential map, i.e., such
that
\be
\Phi(k)\bigl(\exp_p(v)\bigr) = \exp_p\bigl(k\cdot v\bigr)
\ee
for all~$v\in B_r(0_p)\subset T_pM$.  Then~$\Phi(k)$ is
a holomorphic isometry of~$g$ for all~$k\in \bbT$.

Now let~$\d\mu$ be Haar measure on~$\bbT$ and choose any
holomorphic mapping~$\psi:B_r(p)\to T_pM\simeq\bbC^n$
with the property that~$\psi(p)=0$ and~$\psi'(p):T_pM\to T_{0_p}(T_pM)$
is the inverse of the exponential 
mapping~$\exp'_p:T_{0_p}(T_pM)\to T_pM$.  (It may be necessary
to shrink~$r$ to do this.) 

Define a holomorphic mapping~$w:B_r(p)\to T_pM$ by the
averaging formula
\be
w(z) = \int_\bbT k^{-1}\cdot\psi\bigr(\Phi(k)z\bigl)\,\d\mu
\ee
for~$z\in B_r(p)$.
Then~$w(p) = 0_p$ and, by construction,~$w\bigl(\Phi(k)z\bigr)
= k\cdot w(z)$ for all~$z\in B_r(p)$ and all~$k\in\bbT$.   
Moreover, also by construction, $w'(p) = \psi'(p)$.  
In particular, by shrinking~$r$ again, if necessary, 
it can be assumed that~$w$ defines a $\bbT$-equivariant 
holomorphic embedding of~$B_r(p)$ into~$T_pM\simeq\bbC^n$.

In particular, the holomorphic mapping~$w:B_r(p)\to T_pM$
satisfies
\be
w\bigl(\exp_{tY}(z)\bigr) = \exp(tY'_p)\bigl(w(z)\bigr),
\ee
for all real~$t$. Since~$w$ is holomorphic and~$Y$ is 
the imaginary part of the holomorphic vector field~$Z$,
it follows that, for~$z\in B_r(p)$ and~$t$ complex and
of sufficiently small modulus, the identity
\be
w\bigl(\exp_{tZ}(z)\bigr) = \exp(tZ'_p)\bigl(w(z)\bigr)
\ee
holds.  In particular, the coordinate system~$w$ linearizes~$Z$
at~$p$.  
\end{proof}

\begin{remark}[The exponential map]
Of course, the exponential map~$\exp_p:T_pM\to M$ of~$g$ also
intertwines the flow of~$Y'_p$ on~$T_pM$ with the flow of~$Y$
on~$M$, but the exponential map is not generally holomorphic
and so cannot be used to linearize~$Z$ holomorphically.
\end{remark}

\begin{remark}[Complex vs. real flows]
The reader may want to remember that, for a holomorphic
vector field~$Z = X - \iC Y$, the two real vector fields~$X$
and~$Y$ have commuting flows and that, moreover, the identity
\be
\exp_{(a+\iC b)Z} = \exp_{2aX}\circ\exp_{2bY}
\ee
holds. (The factors of~$2$ are neglected in some references.)
\end{remark}

\begin{corollary}
Let~$g$ be a gradient K\"ahler Ricci soliton on~$M$ and
let~$Z$ be its associated holomorphic vector field.
Let~$p\in M$ be a singular point of~$Z$ and let~$\lambda\in\bbR^*$
be a nonzero eigenvalue of~$Z'_p$ of multiplicity~$k\ge1$.
Then there exists a $k$-dimensional complex submanifold~$N_\lambda
\subset M$ that passes through~$p$, to which~$Z$ is everywhere
tangent, and on which~$Y$ is periodic of period~$4\pi/|\lambda|$. \qed
\end{corollary}

\begin{remark}[Nonuniqueness of the $N_\lambda$]
The reader should be careful not to confuse the submanifolds 
$N_\lambda$ with the images under the exponential mapping of the
eigenspaces of~$Z'_p$ acting on~$T_pM$.  Indeed, the~$N_\lambda$
need not be unique.  For example, for the linear vector field
\be
Z = z^1\,\frac{\p\hfill}{\p z^1} 
      + 2z^2\,\frac{\p\hfill}{\p z^2}.
\ee
on~$\bbC^2$, each of the parabolas~$z^2 - c(z^1)^2=0$ for~$c\in\bbC$
is tangent to~$Z$ and the imaginary part of~$Z$ has period~$4\pi$
on all of~$\bbC^2$, so each could be regarded as~$N_1$.  

On the other hand,
the line~$z^1=0$ is the only curve that could be regarded as~$N_2$,
since this is the union of the $2\pi$-periodic points of~$Y$.
\end{remark}

\begin{remark}[Existence at singular points]
\label{rem: singexistexamples}
Example~\ref{ex: products} shows that diagonalizability
with real eigenvalues is sufficient for a linear vector field
to be the linear part of a vector field associated to a 
(locally defined) gradient K\"ahler Ricci soliton.
\end{remark}

\subsubsection{Prescribed eigenvalues}
Let~$\hs = (h_1,\ldots,h_n)\in\bbR^n$ be a nonzero real vector
and define
\be
\Lambda_{\hs} = \left\{\ks\in\bbZ^n\ \vrule\ \ks\cdot\hs=0\right\}
= \bbZ^n\cap \hs^\perp\subset\bbR^n.
\ee
Then~$\Lambda_{\hs}$ is a free abelian group of rank~$n-k$
for some~$1\le k\le n$.  The number~$k$ is the dimension over~$\bbQ$
of the~$\bbQ$-span of the numbers~$h_1,\ldots,h_n$ in~$\bbR$.
Let~$\Lambda^+_\hs\subset\Lambda_\hs$ consist of the~$\ks\in
\Lambda_\hs$ such that~$\ks = (k_1,\ldots,k_n)$ with each
$k_i$ nonnegative.

Consider the linear holomorphic vector field
\be\label{eq: Z_hsdef}
Z_\hs = \sum_{j=1}^n h_j z^j\,\frac{\p\hfill}{\p z^j}
\ee
on~$\bbC^n$.  Let~$Z_\hs = X_\hs - \iC Y_\hs$ be the decomposition 
into real and imaginary parts.  

The closure of the flow of~$Y_\hs$ 
is a connected compact abelian subgroup~$\bbT_\hs\subset\Un(n)$ 
of dimension~$k$. (In fact, in these coordinates, $\bbT_\hs$ 
lies in the diagonal matrices in~$\Un(n)$.)  Note that~$Z_\hs$
and hence~$X_\hs$ are invariant under the action of~$\bbT_\hs$.

\subsubsection{Normalizing volume forms}
In addition to knowing that~$Z$ can be linearized
near a singular point, it will be useful to know
that this can be done in such a way that it simplifies
the coordinate expression for~$\Upsilon$ as well:

\begin{proposition}[Volume normalization at $Z$-singular points]
\label{prop: ZlinUpsconst}
Set~$h = h_1+\cdots+h_n$ and let~$\Upsilon$
be a nonvanishing holomorphic~$n$-form defined 
on an open neighborhood~$U$ of the origin in~$\bbC^n$ 
that satisfies~$\d(Z_\hs\lhk\Upsilon) = h\Upsilon$.

Then there exist~$Z_\hs$-linearizing coordinates~$w = (w^i)$ 
near the origin in~$\bbC^n$ such that, on the domain of
these coordinates~$\Upsilon = \d w^1\w\cdots\w\d w^n$.
\end{proposition}

\begin{proof}
There exists a nonvanishing holomorphic function~$F$ 
on~$U$ that satisfies
\be
\Upsilon = F(z)\,\d z^1\w\cdots\w \d z^n
\ee 
and the function~$F$ must be invariant under the flow of~$Z_\hs$.  
In particular, it follows that~$F$ has a power series expansion
of the form
\be\label{eq: Fseries}
F(z) = c_0 + \sum_{\ks\in\Lambda^+_\hs\setminus\{0\}} c_\ks\, z^\ks
\ee
where~$z^\ks$ is the monomial~$(z^1)^{k_1}\cdots(z^n)^{k_n}$
when~$\ks = (k_1,\ldots,k_n)$ and the $c_\ks$ are constants,
with~$c_0\not=0$ (since, by hypothesis~$F(0)\not=0$).

Now, the series
\be\label{eq: Gseries}
G(z) = c_0 + \sum_{\ks\in\Lambda^+_\hs\setminus\{0\}} 
\frac{c_\ks}{(k_1{+}1)}\, z^\ks
\ee
converges on the same polycylinder that the
series~\eqref{eq: Fseries} does.  The resulting
holomorphic function~$G$ is evidently invariant 
under the flow of~$Z_\hs$ and satisfies 
\be\label{eq: Gdiffeq}
G + z^1\,\frac{\p G}{\p z^1} = F.
\ee
Because~$G$ satisfies~\eqref{eq: Gdiffeq},
the function~$w^1 = z^1G(z)$ satisfies
\be
\d w^1\w\d z^2\w\cdots\w\d z^n
= F(z)\,\d z^1\w\d z^2\w\cdots\w\d z^n.
\ee
Moreover, since~$G$ is~$Z_\hs$-invariant, the function~$w^1$ 
satisfies~$\Lie_{Z_\hs} w^1 = h_1\,w^1$.  

Thus, replacing~$z^1$ by~$w^1$ in the coordinate
chart results in a new $Z_\hs$-linearizing coordinate chart
in which $\Upsilon = \d z^1\w\cdots\w\d z^n$.
\end{proof}

\begin{corollary}[Local normal form near singular points]
\label{cor: singlocnormform}
Let~$Z$ and $\Upsilon$ be a holomorphic vector field and
volume form, respectively on a complex $n$-manifold~$M$.
Let~$p\in M$ be a singular point of~$Z$.  

If there exists a gradient K\"ahler Ricci soliton~$g$ 
with Ricci potential~$f$ on a neighborhood of~$p$ 
whose associated holomorphic vector field 
and volume form are~$Z$ and~$\Upsilon$, respectively,
then there exists an~$\hs\in\bbR^n$ and a $p$-centered
holomorphic chart~$z = (z^i):U\to\bbC^n$ such that,
on~$U$,
\be
Z = h_i\,z^i\,\frac{\p\hfill}{\p z^i}
\qquad\text{and}\qquad
\Upsilon = \d z = \d z^1\w\cdots\w\d z^n.
\ee
\end{corollary}

\begin{proof}
Apply Propositions~\ref{prop: Zlin} 
and~\ref{prop: ZlinUpsconst}.
\end{proof}

\subsubsection{Local solitons near a singular point}
In view of Corollary~\ref{cor: singlocnormform}, 
questions about the local existence and generality 
of gradient K\"ahler Ricci solitons with 
prescribed~$Z$ and~$\Upsilon$ near a singular
point of~$Z$ can be reduced by a holomorphic
change of variables to the study of solitons 
on an open neighborhood of~$0\in\bbC^n$ 
with~$Z=Z_\hs$ for some~$\hs\not=0$
and~$\Upsilon = \d z =\d z^1\w\cdots\w \d z^n$.

\begin{proposition}[Solitons with a prescribed singularity]
\label{prop: singMAforh}
Let~$\phi$ be a strictly pseudo-convex function defined
on a $\bbT_\hs$-invariant, contractible neighborhood 
of~$0\in\bbC^n$ that satisfies
\be\label{eq: ZsingMAeq}
\det\left(\frac{\p^2\phi}{\p z^i\p \zbar^j}\right)
\E^{\frac12\d\phi(X_\hs)} = 1
\ee
and
\be\label{eq: ZsingYinv}
\d\phi(Y_\hs) = 0.
\ee
Then~$\Omega = \frac\iC2\,\p\bar\p\phi$ is the
associated K\"ahler form of a gradient K\"ahler Ricci
soliton with Ricci potential~$f = \frac12\d\phi(X_\hs)$
whose associated holomorphic vector field and volume
form are~$Z_\hs$ and~$\d z^1\w\cdots\w\d z^n$, respectively.

Conversely, if~$g$ is a gradient K\"ahler Ricci soliton
defined on a $\bbT_\hs$-invariant, contractible neighborhood 
of~$0\in\bbC^n$ and~$f$ is a Ricci potential for~$g$ 
that satisfies~$f(0)=0$ such that the associated 
holomorphic vector field and volume
form are~$Z_\hs$ and~$\d z^1\w\cdots\w\d z^n$, respectively,
then~$g$ has a K\"ahler potential~$\phi$ 
that satisfies~\eqref{eq: ZsingMAeq} and~\eqref{eq: ZsingYinv}.
\end{proposition}

\begin{proof}
The first part of the proposition follows by computation,
so nothing more needs to be said.  It remains to establish
the converse statement.

Thus, consider a gradient K\"ahler Ricci soliton~$g$ 
defined on a $\bbT_\hs$-invariant, contractible 
neighborhood~$U$ of~$0\in\bbC^n$ 
with Ricci potential~$f$ satisfying~$f(0)=0$ 
whose associated holomorphic volume form and vector field 
are~$\Upsilon=\d z$ and~$Z_\hs$, respectively. 

The metric~$g$ will necessarily be invariant under~$\bbT_\hs$, 
as will its associated K\"ahler form~$\Omega$.
Since~$\Omega^n = n!\,\iC^{n^2}\,2^{-n}\,\E^{-f}
\,\Upsilon\w\overline{\Upsilon}$,
it follows that~$f$, too, must be invariant under~$\bbT_\hs$.  

On~$U$, there will exist some K\"ahler potential~$\phi$ 
so that~$\Omega = \frac\iC2\p\bar\p\phi$. 
By averaging~$\phi$ over~$\bbT_\hs$, it
can be assumed that~$\phi$ is~$\bbT_\hs$-invariant.
By subtracting a constant, it can be 
assumed that~$\phi(0)=0$.

As has been already noted in~\S\ref{ssec: locpot}, 
the difference~$F = 2f - \d\phi(Z_\hs)$
is a holomorphic function on~$U$. By construction,
$F$ is also necessarily~$\bbT_\hs$-invariant and
vanishes at~$0$. Since~$\phi$ is~$\bbT_\hs$-invariant, 
it follows that $\d\phi(Y_\hs)=0$. 
Thus~$F = 2f - \d\phi(Z_\hs) = 2f - \d\phi(X_\hs)$ 
is real-valued and holomorphic and therefore constant.
Thus,~$F$ vanishes identically, i.e., $f = \frac12\d\phi(X_\hs)$.

Now, however, by construction,~$\phi$ satisfies~\eqref{eq: ZsingMAeq}
and, since~$\phi$ is invariant under the flow of~$Y_\hs$,
it also satisfies~\eqref{eq: ZsingYinv}.
\end{proof}

\begin{remark}[Analyticity in the singular case]
The equation~\eqref{eq: ZsingMAeq} is a $\bbT_\hs$-invariant
real-analytic Monge-Amp\`ere equation whose linearization 
at a strictly pseudo-convex solution~$\phi$ is given by
\be\label{eq: MAlinearizedatphi}
\Delta u + 2\,\Lie_{X_\hs} u = 0
\ee 
where~$\Delta$ is the Laplacian with respect to the metric~$g$
associated to~$\Omega = \frac\iC2\,\p\bar\p\phi$.  Of course,
this is an elliptic equation.

It follows by elliptic regularity that any gradient K\"ahler Ricci 
soliton is real-analytic, even in the neighborhood of singular 
points of~$Z$.  
\end{remark}

\begin{example}[Existence with presribed~$\hs$]
By considering Example~\ref{ex: products}, one sees that,
for any~$\hs$, there is a sufficiently small ball around 
the origin on which there is at least one strictly 
pseudo-convex solution~$\phi$ to~\eqref{eq: ZsingMAeq}.
\end{example}

\subsubsection{A boundary value formulation}
Suppose now that~$\phi$ is a strictly pseudo-convex
solution of~\eqref{eq: ZsingMAeq} defined
on a~$\bbT_\hs$-invariant bounded neighborhood~$D\subset\bbC^n$
of~$0\in\bbC^n$ with smooth boundary~$\p D$.  
Let~$g$ be the corresponding gradient K\"ahler Ricci soliton.

Any solution~$u$ of~\eqref{eq: MAlinearizedatphi}
in~$D$ that vanishes on the boundary will also satisfy
\be
0 = \int_D |\nabla u|^2 + {\ts\frac12}R(g)\,u^2\,\,\d{vol}_g\,,
\ee
as follows by integration by parts using the identities
$\rho = \Lie_{X_h}\Omega$ and~$\d{vol}_g = \frac1{n!}\Omega^n$.

In particular, by shrinking~$D$ if necessary, 
it can be assumed that any solution~$u$ 
to~\eqref{eq: MAlinearizedatphi} in~$D$ that vanishes 
on~$\p D$ must vanish on~$D$.

It then follows, by the implicit function theorem,
that any~$\bbT_\hs$-invariant function~$\psi$ on~$\p D$
that is sufficiently close (in the appropriate norm) 
to~$\phi$ on~$\p D$ is the boundary value of a unique
pseudo-convex solution~$\tilde\phi$ of~\eqref{eq: ZsingMAeq} 
that is near~$\phi$ on~$D$. The uniqueness then implies
that~$\tilde\phi$ must also be~$\bbT_\hs$-invariant
and so must, in particular, satisfy~\eqref{eq: ZsingYinv}.

Note that the metric~$g$ does not always uniquely 
determine~$\phi$ by the construction given
in Proposition~\ref{prop: singMAforh} since one can add 
to~$\phi$ the real part of any $\bbT_\hs$-invariant 
holomorphic function that vanishes at~$0\in\bbC^n$.  
(Depending on~$\hs$, there may or may not be any 
nonconstant $\bbT_\hs$-invariant holomorphic functions
on a neighborhood of~$0\in\bbC^n$.)  However, this 
ambiguity is relatively small.

Thus, local gradient K\"ahler Ricci solitons near~$0\in\bbC^n$ 
with prescribed holomorphic data~$(Z,\Upsilon) = (Z_\hs,\d z)$
do exist and have a `degree of generality' that
depends on the number~$k$.  The most constraints
appear when~$k$ reaches its maximum value~$n$ 
and the least when~$k$ reaches its minimum value~$1$.

\section{Poincar\'e coordinates in the positive case}
\label{sec: poincare}

Throughout this section,~$M$ will be a noncompact, 
simply connected complex manifold and~$g$ will be 
a complete gradient K\"ahler Ricci soliton with 
postive Ricci curvature.  Moreover, it will be assumed 
that the scalar curvature~$R(g)$ has at least one 
critical point. 

\subsection{First consequences}
Cao and Hamilton~\cite[Proposition~4.2]{MR01g:53117} prove 
the following useful result:

\begin{lemma}\label{lem: uniqueRcrit}
The scalar curvature~$R(g)$ has only one critical point
and it is both a local maximum and the unique critical
point of~$f$, which is a strictly convex proper function on~$M$.
\end{lemma}

\begin{proof}
Since~$R(g) + 2|Z|^2 = 2h$ by Proposition~\ref{prop: RnormZ2reln}, 
the function~$R(g)\ge0$ is bounded by the constant~$h$ 
and any critical point of~$R(g)$ is a critical point 
of~$|Z|^2 = \frac12|\nabla f|^2$.  

On the other hand, since~$\nabla^2f=\Ric(g)$,
which is positive definite, the formula
\be
\d\left({\ts\frac12}|\nabla f|^2\right)(\nabla f)
= \nabla^2f(\nabla f,\nabla f) = \Ric(g)(\nabla f,\nabla f)
\ee
shows that~$\frac12|\nabla f|^2$ cannot have any critical point
away from where~$\nabla f = 0$.  Moreover, any point~$p$
where~$\nabla f$ vanishes satisfies $R(g)(p) = 2h$, 
which is the maximum possible value of~$R(g)$.

Since~$\nabla^2f=\Ric(g)$ is positive definite,
the function~$f$ is locally strictly convex.  
Since~$g$ is complete,~$f$ can have at most one critical
point, i.e., point where~$\nabla f = 0$, and it must be a 
nondegenerate minimum of~$f$.  

By hypothesis, there does exist a (unique) critical point 
of~$f$; call it~$p$.  By adding a constant to~$f$ it can 
be assumed that~$f(p)=0$.  It remains to show that~$f$ is
proper, i.e., that~$f^{-1}\bigl([a,b]\bigr)\subset M$
is compact for any closed interval~$[a,b]\subset\bbR$.

Since~$R(g) + 2|Z|^2 = 2h$ and since~$R(g)>0$, 
it follows that~$|Z|\le\sqrt{h}$, so that~$Z$ has
bounded length.  In particular, writing
\be
Z=X-\iC Y = {\ts\frac12}\bigl(\nabla f - \iC J(\nabla f)\bigr) 
\ee
one has~$|X|^2 = |Y|^2 = \frac12|Z|^2 < \frac12h$, 
so $X$ and~$Y$ have bounded lengths as well. Since~$g$
is complete, their flows are globally defined on~$M$.

Let~$\gamma:\bbR\to M$ be any nonconstant integral curve 
of~$\nabla f$, i.e.,~$\gamma'(t) 
= \nabla f\bigl(\gamma(t)\bigr)\not=0$
for all~$t\in\bbR$.  Consider the function~$\phi(t) 
= f\bigl(\gamma(t)\bigr)$.  Straightforward computation yields
$\phi'(t) = \bigl|\nabla f\bigl(\gamma(t)\bigr)\bigr|^2>0$
and
\be\label{eq: phi2ndpos}
\phi''(t) = 2\Ric(g)\bigl(\nabla f\bigl(\gamma(t)\bigr),\nabla f\bigl(\gamma(t)\bigr)\bigr) >0,
\ee
so~$\phi:\bbR\to\bbR$ is strictly convex and increasing.  
It follows that~$\phi$ increases without bound along~$\gamma$.

Since~$\nabla^2f$ is positive definite, the
critical point~$p$ is a source singularity of the 
vector field~$\nabla f$.  Let~$U\subset M$ be the open set 
that consists of~$p$ and all of the points~$q$ in~$M$ whose
$\alpha$-limit point under~$\nabla f$ is equal to~$p$.  
Since~$f$ strictly increases without bound on each integral 
curve of~$\nabla f$, it follows that~$f$ maps each
integral curve of~$\nabla f$ that lies in~$U$ diffeomorphically
onto~$(0,\infty)$.  Moreover, for each~$c>0$, the 
set~$f^{-1}(c)\cap U$ is compact and diffeomorphic to~$S^{2n-1}$.
Indeed, $f:U\to[0,\infty)$ is proper.

Now suppose that~$U\not=M$.  Then, by the connectedness of~$M$, 
there exists a point~$q\in M\setminus U$ that is not in the interior
of~$M\setminus U$, i.e., a point~$q\not\in U$ such that there
exists a sequence~$q_i\in U$ that converges to~$q$.
This implies, in particular, that~$f(q_i)\ge0$ converges to~$f(q)=c$.
Thus, $c\ge 0$ and, for~$i$ sufficiently large, 
$q_i$ must lie in~$f^{-1}\bigl([0,c{+}1]\bigr)\cap U$, 
which has been shown to be compact 
and must therefore contain its limit points. 
Thus~$q$ lies in~$f^{-1}\bigl([0,c{+}1]\bigr)\cap U$,
although, by construction,~$q\not\in U$.  Thus, $U=M$
and~$f$ is proper, as claimed.
\end{proof}

\begin{remark}[$M$ is Stein]
As Cao and Hamilton remark, since~$\rho=\iC\p\bar\p f$
is the Ricci form of~$g$, which is positive, the proof shows 
that~$f$ is a strictly plurisubharmonic proper exhaustion 
function on~$M$.  This implies that~$M$ is Stein and, 
as Cao points out in~\cite[Proposition 3.2]{math.DG/0302112},
that~$M$ is diffeomorphic to~$\bbR^{2n}$.

However, as will be seen in Theorem~\ref{thm: poincarecoords},
one has the stronger result that 
$M$ is biholomorphic to~$\bbC^n$.
\end{remark}

The following result, also known to Cao and Hamilton,%
\footnote{H.-D. Cao, personal communication, 2 June 2004}
gives constraints on the rate of growth of the Ricci potential.

\begin{lemma}[Growth of~$f$]\label{lem: fgrowth}
Let~$p$ be the critical point of~$R(g)$ and let~$f$
be the Ricci potential, normalized so that~$f(p)=0$.
There exist positive constants~$c_1$ and~$c_2$ such
that, for all~$x\in M$,
\be
\sqrt{1+\bigl(c_1\,d(x,p)\bigr)^2}-1 \le f(x) \le c_2\,d(x,p).
\ee
\end{lemma}

\begin{proof}
Since~$g$ is complete, there exists a geodesic
joining~$p$ to~$x$ whose length is~$d(p,x)$.
Let~$\alpha:\bbR\to M$ be such a unit speed geodesic
with~$\alpha(0)=p$ and~$\alpha(s)=x$ such that
$d(p,x) = s$.  

Consider the function~$\phi(t) = f\bigl(\alpha(t)\bigr)$.
By the Chain Rule, and the fact that~$\alpha$ has unit speed,
\be
\phi'(t) = \nabla f\bigl(\alpha(t)\bigr)\cdot\alpha'(t)
 \le \bigl|\nabla f\bigl(\alpha(t)\bigr)\bigr| \le \sqrt{2h}.
\ee
Since~$\phi(0)=0$, it follows 
that~$f(x) = f\bigl(\alpha(s)\bigr) = \phi(s)\le \sqrt{2h}\,s$.
Thus, one can take~$c_2 = \sqrt{2h}$.  

For the other inequality, note that, again, by the Chain Rule,
\be
\phi''(t) 
= \nabla^2 f\bigl(\alpha(t)\bigr)\bigl(\alpha'(t),\alpha'(t)\bigr)
= \Ric(g)\bigl(\alpha(t)\bigr)\bigl(\alpha'(t),\alpha'(t)\bigr)
\ee
and the right hand side of this equation is positive since~$\Ric(g)$
is positive.  Moreover, if~$\lambda_{\textrm{min}}(g)>0$
denotes the minimum eigenvalue of~$\Ric(g)$, which is
a positive continuous function on~$M$, it follows that
\be
\phi''(t) \ge \lambda_{\textrm{min}}(g)\bigl(\alpha(t)\bigr) > 0.
\ee
In particular,~$\phi$ is a convex function on~$\bbR$.

Let~$r_0>0$ be sufficiently small that it is below the 
injectivity radius of~$g$ at~$p$ and sufficiently small
that~$\lambda_{\textrm{min}}(g)(y)
\ge\frac12\,\lambda_{\textrm{min}}(g)(p)$ for all~$y$
lying within~$B_{r_0}(p)$.  Let~$a = \frac12\,\lambda_{\textrm{min}}(g)(p)>0$.

Then~$\phi''(t)\ge a$ for~$|t|\le r_0$ while 
$\phi''(t)>0$ for~$|t|\ge r_0$. Because~$\phi(0)=\phi'(0)=0$,
it follows that $\phi(t)\ge A(t)$ for all~$t\in\bbR$ 
where
\be
A(t) = 
\begin{cases}
 {\ts\frac12}at^2 & \text{for $|t|\le r_0$;}\\
 ar_0|t|- {\ts\frac12}a\,{r_0}^2 & \text{for $|t|\ge r_0$.}
\end{cases}
\ee
Since there exists a positive constant~$c_1$ such that
$A(t)\ge \sqrt{1+(c_1t)^2}-1$, the desired lower bound
follows.
\end{proof}

\begin{remark}[An alternative growth formulation]
Another formulation of~Lemma~\ref{lem: fgrowth} is 
that the function~$c:M\setminus\{p\}\to\bbR$ 
defined by
\be
c(x) = \frac{\sqrt{f(x)\bigl(f(x)+2\bigr)}}{d(x,p)}>0
\ee
is bounded above and has a positive lower bound.
\end{remark}

The bounds of~Lemma~\ref{lem: fgrowth} can be simplified 
somewhat if one stays sufficiently far from~$p$:

\begin{corollary}\label{cor: flingrowth}
For every~$r>0$, there exist positive constants~$c_1$
and~$c_2$ such that, for all~$x$ outside the ball of 
radius~$r$, one has
\be
c_1\,d(x,p) \le f(x) \le c_2\,d(x,p).
\ee \qed
\end{corollary}

\begin{remark}[The growth rate of~$f$]
For any vector~$v\in TM$, one has
\be
\Ric(g)(v,v) \le \lambda_{\textrm{max}}(g)\,|v|^2
\ee
where~$\lambda_{\textrm{max}}(g):M\to\bbR$ is the
maximum eigenvalue function for $\Ric(g)$.
Since~$g$ is K\"ahler, the eigenvalues of~$\Ric(g)$
occur in pairs and, since~$\Ric(g)>0$, it follows that
$\lambda_{\textrm{max}}(g)\le \frac12 R(g)$. In particular,
by Proposition~\ref{prop: RnormZ2reln},
one has the more explicit inequality
\be\label{eq: Ricbound}
\Ric(g)(v,v) \le {\ts\frac12}R(g)\,|v|^2
\le {\ts\frac12}\bigl(2h - |\nabla f|^2\bigr)\,|v|^2.
\ee

Now let~$\gamma:(0,\infty)\to M$ be the arclength
parametrization of a nonconstant integral curve of~$\nabla f$,
such that~$p$ is the limit of~$\gamma(s)$ as~$s\to 0^+$.
Thus,~$\bigr|\nabla f\bigl(\gamma(s)\bigr)\bigr|\,\gamma'(s) 
= \nabla f\bigl(\gamma(s)\bigr)$ for all~$s>0$.

Let~$\phi(s) = f\bigl(\gamma(s)\bigr)$. One then
computes via the Chain Rule that 
\be
\phi'(s) = \bigl|\nabla f\bigl(\gamma(s)\bigr)\bigr|
\le \sqrt{2h}.
\ee
and hence that
\be\label{eq: phiprimeprime}
\phi''(s)
= \Ric(g)\left(\frac{\nabla f\bigl(\gamma(s)\bigr)}
             {\left|\nabla f\bigl(\gamma(s)\bigr)\right|},
\frac{\nabla f\bigl(\gamma(s)\bigr)}
             {\left|\nabla f\bigl(\gamma(s)\bigr)\right|}\right).
\ee
By the positivity of~$\Ric(g)$ and~\eqref{eq: Ricbound}, 
this implies
\be
0< \phi''(s)  
\le {\ts\frac12}\left(2h - \left(\phi'(s)\right)^2\right).
\ee
Moreover, it is clear that, as~$s\to0^+$, the quantity
on the right hand side of~\eqref{eq: phiprimeprime}
has~$\lambda_{\textrm{min}}(g)(0)>0$ as a lower bound for
its infimum limit.  Thus, the infimum limit of~$\phi''(s)$
as~$s\to0^+$ is positive.

From these relations, several conclusions can be drawn.
The function~$\phi$ is increasing and strictly convex up
on~$(0,\infty)$.  On the other hand, since~$\phi'$ is 
bounded above, it follows that~$\phi$ grows at most
linearly. Moreover, there must be a sequence
of distances~$s_k\to\infty$ such that~$\phi''(s_k)\to 0$.
Since, by~\eqref{eq: phiprimeprime}
\be
\phi''(s_k) \ge \lambda_{\textrm{min}}(g)\bigl(\gamma(s_k)\bigr),
\ee 
it follows that $\lambda_{\textrm{min}}(g)\bigl(\gamma(s_k)\bigr)
\to0$ as~$k\to\infty$.
\end{remark}

\subsection{Poincar\'e coordinates}
Let~$\Upsilon$ be the associated holomorphic volume
form on~$M$, normalized so that $\Upsilon$ has
unit size at~$p$.  This determines~$\Upsilon$ up to
a complex multiple of modulus~$1$.  Let~$Z$ be the
associated holomorphic vector field.

Since~$Z$ vanishes at~$p$, the eigenvalues of~$Z'_p$
are the eigenvalues of the Ricci tensor at~$p$, which
are real and positive, say~$h_1,\ldots,h_n>0$.
Set~$h = h_1+\cdots+h_n>0$, as usual.  

\begin{theorem}[Poincar\'e coordinates]\label{thm: poincarecoords}
There exists a global special coordinate system~$z:M\to\bbC^n$
that linearizes~$Z$.  In particular, $M$ is biholomorphic to~$\bbC^n$.
\end{theorem}

\begin{proof}
By Proposition~\ref{prop: Zlin}, there exists a small
open ball~$U$ about~$p$ on which there exist $p$-centered 
holomorphic coordinates~$w = (w^i):U\to\bbC^n$ that 
linearize~$Z$. By shrinking~$U$ if necessary, it can be
assumed that~$U = f^{-1}\bigl([0,\varepsilon)\bigr)$ for
some small~$\varepsilon>0$.  Note that, since the~$w^i$
linearize~$Z$, the identity
\be\label{eq: exptZinw}
w^i\bigl(\exp_{tZ}(q)\bigr) = \E^{h_it}\,w^i(q)
\ee
holds for all~$q\in U$ and all~$t\in\bbC$ in the 
connected domain containing~$0\in\bbC$ for 
which~$\exp_{tZ}(q)$ lies in~$U$.  In particular,
this implies that 
\be\label{eq: exptXinw}
w^i\bigl(\exp_{2tX}(q)\bigr) = \E^{h_it}\,w^i(q)
\ee
for all~$q\in U$ and all~$t\in\bbR$ in the interval 
containing~$0\in\bbR$ for which~$\exp_{2tX}(q)$ 
lies in~$U$.

Now, for~$q\in M$ distinct from~$p$, write~$q = \exp_{2t'X}(q')$ 
for some~$q'\in U$ and~$t'\in\bbR$. Define
\be\label{eq: defzcoords}
z^i(q) = \E^{h_it'}w^i(q').
\ee
If~$\exp_{2t'X}(q')=\exp_{2t''X}(q'')$ for 
some~$q''\in U$ and~$t''\in\bbR$, then one
sees from~\eqref{eq: exptXinw} that $\E^{h_it''}w^i(q'')
= \E^{h_it'}w^i(q')$, so~$z^i(q)$ is well-defined. 

Since the flow of~$X$ is holomorphic and~$w^i$ is
holomorphic on~$U$, the function~$z^i:M\to\bbC$ 
is also holomorphic.  Moreover, by construction, 
\be\label{eq: exptXinz}
z^i\bigl(\exp_{2tX}(q)\bigr) = \E^{h_it}\,z^i(q)
\ee
for all~$q\in M$, which implies that
\be\label{eq: exptZinz}
z^i\bigl(\exp_{tZ}(q)\bigr) = \E^{h_it}\,z^i(q).
\ee
In particular, the Lie derivative of~$z^i$ by~$Z$ is $h_i z^i$.

The fact that the mapping~$z = (z^i):M\to\bbC^n$ is one-to-one
and onto now follows immediately since, as was observed in
the proof of~Lemma~\ref{lem: uniqueRcrit}, the gradient flow lines
of~$\nabla f = 2X$ all have~$p$ as $\alpha$-limit point and
the flow of~$\nabla f$ exists for all time. 

Finally, in these coordinates~$\Upsilon = F(z)\,\d z^1\w\cdots\w\d z^n$
for some nonvanishing entire holomorphic function~$F$ on~$\bbC^n$.
However, since~$\d(Z\lhk\Upsilon) = h\,\Upsilon$, it follows
immediately that~$\d F(Z) = 0$.  Since all of the eigenvalues
of~$Z'_p$ are positive, this is only possible if~$F$ is a constant
function.  By scaling one of the~$z^i$ by a constant, 
it can be arranged that~$F\equiv1$.

Thus, the resulting global coordinate system~$(M,z)$ 
is special and linearizes~$Z$, as desired.
\end{proof}

\begin{remark}[Previous results]
Chau and Tam~\cite[Theorem 1.1]{math.DG/0310198} 
proved that~$M$ is biholomorphic to~$\bbC^n$ under the 
additional hypothesis that all the eigenvalues~$h_i$ 
are equal.  In a very recent posting to the 
arXiv~\cite{math.DG/0404449}, they prove a result that
implies that~$M$ is biholomorphic to~$\bbC^n$ under the
hypotheses of~Theorem~\ref{thm: poincarecoords}.  However,
their result does not provide $Z$-linearizing coordinates,
which is the main purpose of~Theorem~\ref{thm: poincarecoords}.
\end{remark}

\subsection{Coordinate ambiguities}
The reader may find it surprising that any
local $Z$-linearizing coordinates~$z^i$ defined
on a neighborhood of the $Z$-singular point~$p$ extend 
to global coordinates on~$\bbC^n$ that are special for
any gradient K\"aher-Ricci soliton defined on~$\bbC^n$
with positive Ricci curvature whose associated 
holomorphic vector field is~$Z$.  

This is perhaps made less surprising by the following
result:

\begin{proposition}
Let~$\hs = (h_1,\ldots,h_n)\in \bbR^n$
be a vector with~$h_i>0$ for~$1\le i\le n$.
Consider the vector field
\be\label{eq: Z_hs}
Z_\hs = h_i\,z^i\,\frac{\p\hfill}{\p z^i}
\ee
on~$\bbC^n$.  Then the set~$G_\hs$ of 
biholomorphisms~$\psi:\bbC^n\to\bbC^n$ that preserve~$Z_\hs$ 
is a complex Lie group of dimension~$d_\hs$ where~$d_\hs\ge n$ 
is the number of vectors~$\ks = (k_1,\ldots,k_n)\in\bbZ^n$ 
that satisfy~$k_i\ge0$ and~$\ks\cdot\hs\in\{h_1,\ldots,h_n\}$.

Moreover, if~$U\subset\bbC^n$ is any connected open
neighborhood of~$0\in\bbC^n$, 
then any locally defined biholomorphism~$\psi:(U,0)\to(\bbC^n,0)$
that preserves~$Z_\hs$ is the restriction to~$U$ of an
element of~$G_\hs$.
\end{proposition}

\begin{proof}
Let~$U\subset\bbC^n$ be an open neighborhood of~$0$
and let~$\psi = \bigl(w^i(z)\bigr):U\to\bbC^n$ be a
local biholomorphism that preserves~$Z$.  Since~$Z$
has only one singular point, namely~$0\in\bbC^n$,
it follows that~$\psi(0)=0$.  Moreover, by construction,
the functions~$w^i$ must satisfy~$\d w^i(Z) = h_i w^i$.
It follows that each~$w^i$ has a power series expansion
about~$0\in\bbC^n$ of the form
\be\label{eq: G_hdefined}
w^i = \sum_{\{\ks\ge0\ \vrule\ \ks\cdot\hs = h_i\}} c^i_\ks\,z^\ks.
\ee
Since the right hand side has only a finite number of
terms, this power series is a polynomial and hence
globally defined on~$\bbC^n$.  It remains to see that
it is invertible.

Consider the $n$-form~$\d w = \d w^1\w\cdots\w\d w^n$.
By the above analysis~$\d w = F(z)\,\d z$ for some 
polynomial~$F(z)$. By hypothesis, $\psi$ is
a local biholomorphism, so~$F(0)\not=0$.  Since~$\Lie_Z\d w
= (h_1+\cdots+h_n)\,\d w$ by construction, it follows
that~$\d F(Z) = 0$, i.e., that~$F$ is~$Z$-invariant.
This implies that~$F$ is constant and hence nowhere
vanishing.

Now, by hypothesis,~$\psi$ is locally invertible,
with, say, a local inverse~$\psi^{-1}:(V,0)\to(\bbC^n,0)$.
However, by construction,~$\psi^{-1}$ preserves~$Z$,
so, by the argument given above, $\psi^{-1}$ is
also a polynomial mapping and hence extends to 
a global polynomial mapping~$\pi:(\bbC^n,0)\to(\bbC^n,0)$.
Since~$\psi\circ\pi:(\bbC^n,0)\to(\bbC^n,0)$ is a
polynomial mapping that is the identity on some 
neighborhood of~$0$, it must be the identity everywhere
on~$\bbC^n$. In particular,~$\pi$ is the
global inverse of~$\psi$ extended to~$\bbC^n$, 
which is now revealed to be an element of~$G_\hs$,
which is what needed to be shown.

Finally, it is clear that, for any~$i$ and any choice 
of~constants~$c^i_\ks\in\bbC$ for~$(i,\ks)$ such that
$\ks\in\bbZ^n$ satisfies~$k_j\ge0$ for~$1\le j\le n$
and~$\ks\cdot\hs = h_i$, the formula~\eqref{eq: G_hdefined}
defines a polynomial~$w^i$ that satisfy~$\Lie_Z w^i = h_i w^i$.

Moreover, for any choice of~$d_\hs$ constants~$c = (c^i_\ks)$
where~$(i,\ks)$ satisfies $\ks\in\bbZ^n$ with~$k_j\ge0$ 
for~$1\le j\le n$ and~$\ks\cdot\hs = h_i$, the corresponding
collection of functions~$w^i$ satisfies
\be
\d w^1\w\cdots\w \d w^n 
= F(c^i_\ks)\,\d z^1\w\cdots\w \d z^n. 
\ee
where~$F$ is a polynomial of degree~$n$ in the~$d$ parameters
$c^i_\ks\in\bbC$.  

As long as~$F(c^i_\ks)\not=0$, the polynomial mapping~$\psi_c = (w^i)$ 
is a local (and therefore global) biholomorphism of~$\bbC^n$ 
that preserves~$Z$ and hence lies in~$G_\hs$.  
Thus, the~$c^i_\ks$ define global holomorphic coordinates 
on~$G_\hs$ that embed it into~$\bbC^{d_\hs}$ as an open set.
\end{proof}

\begin{remark}[The structure of~$G_\hs$]
If~$\mu_1,\ldots,\mu_k\ge 1$ are the multiplicities
of the eigenvalues~$(h_1,\ldots,h_n)$, then~$G_\hs$ 
is the semi-direct product 
of a reductive subgroup isomorphic 
to~$\GL(\mu_1,\bbC)\times\cdots\GL(\mu_k,\bbC)$
with a nilpotent subgroup biholomorphic to~$\bbC^\mu$
where~$\mu = d_\hs - {\mu_1}^2 - \cdots - {\mu_k}^2$.

When~$n=1$, one has~$G_\hs\simeq\bbC^*=\GL(1,\bbC)$.
When~$n=2$, one has either 
\begin{enumerate}
\item $d_\hs =2$ if~$\hs = (h_1,h_2)$ 
with neither~$h_1/h_2$ nor $h_2/h_1$ an integer
(in which case~$G_\hs=\bbC^*\times\bbC^*$);
\item $d_\hs =3$ if~$\hs = (h_1,h_2)$ 
with either~$h_1/h_2$ or $h_2/h_1$ an integer greater than~$1$; or
\item $d_\hs =4$ if~$\hs = (h,h)$ 
(in which case~$G_\hs = \GL(2,\bbC)$).
\end{enumerate}

When~$n>2$, there is no upper bound for~$d_\hs$ that depends only on~$n$.
For example, when~$n=3$, one has~$d_{(1,1,k)} = k+6$ for
any integer~$k>1$.
\end{remark}

\subsection{Global consequences}
Throughout this section, $g$ will be a complete
gradient K\"ahler Ricci soliton on~$\bbC^n$
with positive Ricci curvature whose
associated vector field~$Z$ is given by~\eqref{eq: Z_hs}
where~$\hs = (h_1,\ldots,h_n)$ and
\be
0< h_1 \le h_2 \le \cdots \le h_n\,.
\ee

The compact abelian group~$\bbT_\hs\subset\Un(n)$ will
denote the closure of the orbit of~$Y$, the imaginary
part of~$Z$.

The existence of global linearizing coordinates for
a gradient K\"ahler Ricci soliton gives elementary
proofs and/or improvements of several known results.

\subsubsection{Periodic orbits}
The first result sharpens Theorem~1.1 of the 
article~\cite{MR01g:53117} of Cao and Hamilton.

\begin{proposition}[Periodic orbits of~$\nabla f$]
\label{prop: perorbnablaf}
For all~$c>0$, the flow of~$J(\nabla f)$
preserves the {\upshape(}smooth{\upshape)}
level set~$f^{-1}(c)\subset M$
and has at least~$n$ periodic orbits on~$M$.
\end{proposition}

\begin{proof}
Since $Z = \frac12\bigl(\nabla f - \iC\,J(\nabla f)\bigr)$,
and since~$h_i>0$ for $1\le i\le n$, it follows that
$J(\nabla f)$  is periodic of period~$2\pi/h_i$ on
the~$z^i$-axis.  Moreover, since~$f$ increases without
bound as~$|z^i|\to\infty$, this axis meets each level
set~$f^{-1}(c)$ for~$c>0$ in a circle.  Thus, there
are at least $n$ distinct periodic orbits of~$J(\nabla f)$
within each such level set.
\end{proof}

\subsubsection{An invariant potential}
As has already been seen, the metric $g$ is invariant 
under~$\bbT_\hs$.  It turns out that one can canonically
choose a K\"ahler potential for~$g$:

\begin{proposition}[Canonical potentials]\label{prop: canpot}
There is a unique $\bbT_\hs$-invariant 
K\"ahler potential~$\phi:\bbC^n\to\bbR$
satisfying~$\Omega = \frac\iC2\,\p\bar\p\phi$ 
and~$\phi(0)=0$.
\end{proposition}  

\begin{proof}
Since~$M=\bbC^n$, there exists at least one K\"ahler 
potential~$\phi$ for~$g$, i.e., such that~$\Omega =
\frac\iC2\,\p\bar\p\phi$.  Since~$\bbT_\hs$ is compact,
by averaging~$\phi$ over~$\bbT_\hs$, one can assume that~$\phi$
is~$\bbT_\hs$ invariant and by adding a constant, 
one can assume that~$\phi(0)=0$.  

If~$\tilde\phi$ were also~$\bbT_\hs$-invariant and
satisfied~$\Omega =\frac\iC2\,\p\bar\p\tilde\phi$,
then the difference~$\tilde\phi - \phi$ would be
the real part of a~$\bbT_\hs$-invariant holomorphic
function~$H$.  In particular~$H$ would be invariant
under the flow of~$Y$ and hence of~$Z$.  However, as
has already been seen, the only holomorphic functions on~$\bbC^n$
that are invariant under the flow of~$Z$ are the constants.
Thus~$\tilde\phi-\phi$ is constant.  The normalization~$\phi(0)=0$
then guarantees the uniqueness of~$\phi$.
\end{proof}

\subsubsection{Normalized linearizing coordinates}
The ambiguity in the linearizing
coordinates for the vector field~$Z$
represented by the group~$G_\hs$ 
can be used to simplify the potential for~$g$.

\begin{theorem}[Normalized coordinates]\label{thm: normedcoords}
Let~$\phi$ be the unique~$\bbT_\hs$-invariant 
K\"ahler potential for~$g$, normalized so that~$\phi(0)=0$.
Then there exists an element~$\Psi\in G_\hs$,
unique up to composition with an
element of the compact group~$\Un(n)\cap G_\hs$, 
such that
\be
\Psi^*(\phi) = |z^1|^2 + \cdots + |z^n|^2 
+ E_{\bi\bj kl}(z) \zbar^i\zbar^j z^k z^l
\ee
for some real-analytic functions~$E_{\bi\bj kl}=E_{\bj\bi kl}
=E_{\bi\bj lk}=\overline{E_{\bk\bl ij}}$ 
defined near~$0\in\bbC^n$.
\end{theorem}

\begin{proof}
Let~$f$ be the Ricci potential for~$g$, normalized
so that~$f(0)=0$.  Since~$f$ is~$\bbT_\hs$-invariant
and since, by~\eqref{eq: 2f-pphiZhol}, the difference~$2f-\d\phi(Z)$ 
is holomorphic and~$\bbT_\hs$-invariant, it follows
by the same argument as above that~$2f-\d\phi(Z)$ is
constant and hence vanishes identically.  Thus
\be
\d\phi(Z) = \d\phi(X) = 2 f.
\ee

Because~$\phi$ and~$f$ are real-analytic they have convergent
power series expansions near~$0\in\bbC^n$.  Since~$f(0)=0$
and~$f$ has a critical point at~$0$, it has an expansion
of the form
\be
f = {\ts\frac12}f_{ij}\,z^iz^j+ f_{i\bj}\,z^i\zbar^j
   +{\ts\frac12}\overline{f_{ij}}\,\zbar^i\zbar^j + O(|z|^3).
\ee
where~$f_{ij}=f_{ji}$ and~$f_{i\bj} = \overline{f_{j\bi}}$.
Because of the positivity of the~$h_i$ and the invariance
of~$f$ under the flow of~$Y$, it follows that~$f_{ij}=0$
and~$(h_i{-}h_j)f_{i\bj}=0$ for all~$i$ and~$j$.  
Moreover, since~$f$ is strictly convex up at the
origin, the Hermitian form~$f_{i\bj}\,z^i\zbar^j$
is positive definite.

Thus, by making a linear change of variables that preserves~$Z$
(i.e., by applying a transformation 
in~$\GL(n,\bbC)\cap G_\hs$), 
it can be arranged that
\be
f =  {\ts\frac12}h_1\,|z^1|^2 +\cdots 
      + {\ts\frac12}h_n\,|z^n|^2 + O(|z|^3).
\ee

Next, consider the part of~$f$ that is pure in~$z$ or~$\zbar$,
i.e., consider the expansion
\be
f =  {\ts\frac12}h_1\,|z^1|^2 +\cdots 
      + {\ts\frac12}h_n\,|z^n|^2 
+ \sum_{\ks\ge0,\,|\ks|\ge3} 
\left(f_\ks\, z^\ks + \overline{f_\ks}\,\zbar^\ks\right)
 + f_{i\bj}(z)z^i\zbar^j.
\ee
where~$f_\ks\in\bbC$ and~$f_{i\bj} = \overline{f_{j\bi}}$
vanishes at~$z=0$.  The invariance of~$f$ under the flow
of~$Y$ implies that~$f_\ks = 0$ for all~$\ks$, so these
`pure' terms do not appear after all.

Finally, consider the part of the remainder that is linear
in the variables~$\zbar^i$ or~$z^i$ and vanishes at~$z=0$
to order at least~$3$, i.e., write
\be
f =  {\ts\frac12}h_k\,|z^k|^2 
       + Q^i(z)\,\zbar^i + \overline{Q^i(z)}\,z^i
+ f_{\bi\bj kl}(z) \zbar^i\zbar^j z^k z^l,
\ee
where~$Q^i(z)$ is a holomorphic function of~$z$ that
vanishes to order at least~$2$ at~$z=0$
and~$f_{\bi\bj kl}=f_{\bj\bi kl}
=f_{\bi\bj lk}=\overline{f_{\bk\bl ij}}$.

Again, the fact that~$f$ is invariant under the flow of~$Y$
implies that~$Q^i$ must satisfy~$\Lie_ZQ^i = h_i Q^i$, i.e.,
that~$Q^i$ has an expansion of the form
\be\label{eq: Qiexp}
Q^i(z) = \sum_{\{\ks\ge0\ \vrule\ \ks\cdot\hs = h_i\}}
c^i_\ks z^\ks
\ee
with~$c^i_\ks=0$ unless~$|\ks|=k_1+\cdots+k_n>1$. In particular,
this implies that~$Q^i$ is a polynomial in~$z$ since the
right hand side of~\eqref{eq: Qiexp} can contain only
finitely many terms.  Now consider the change of variables
defined by
\be
w^i = z^i + \frac2{h_i}\,Q^i(z)
\ee
This transformation belongs to~$G_\hs$ by definition and satisfies
\be
f =  {\ts\frac12}h_k\,|w^k|^2 
   + f^*_{\bi\bj kl}(w)\,\wbar^i\wbar^j w^k w^l,
\ee
for some functions~$f^*_{\bi\bj kl}$ with the same symmetry
and reality properties as the corresponding~$f_{\bi\bj kl}$.

Since~$\Lie_X\phi = 2f$ and~$\phi(0)=0$, it follows that~$\phi$
has a power series expansion
\be
\phi = |w^k|^2 + E_{\bi\bj kl}(w)\,\wbar^i\wbar^j w^k w^l,
\ee
as desired.  The uniqueness of the transformation~$\Psi = (w^i)$
up to composition with an element of~$\Un(n)\cap G_\hs$ 
is now evident.
\end{proof}

\subsubsection{Totally geodesic submanifolds}
Since the fixed locus of an isometry of~$g$ must
be totally geodesic, one has the following result:

\begin{proposition}[Geodesic subspaces]\label{prop: geodsubspaces}
If~$h_i$ has multiplicity~$\mu_i>0$ and has the property that, 
for all~$k$, $h_k\not=mh_i$ for any integer~$m>1$, then
the $\mu_i$-plane in~$\bbC^n$ defined by~$z^j=0$ when~$h_j\not=h_i$
is totally geodesic.

More generally, if~$Y$ has a periodic point~$q$ with
period~$T>0$, then the union of the $T$-periodic points
is a nontrivial totally geodesic linear subspace of~$\bbC^n$
generated by the $z^i$-axis lines for which~$h_i$ is 
an integer multiple of~$4\pi/T$. \qed
\end{proposition}

\begin{remark}[Geodesic axes]\label{rem: geodaxes}
The reader might wonder whether or not the hypothesis
of~$h_i$ having no `supermultiples' is necessary in order
for the $h_i$-eigenspace of~$Z_\hs$ in~$\bbC^n$ to be
totally geodesic.  

The answer is clearly `yes' in general $Z$-linearizing
coordinates:  For example, if~$n=2$ and~$\hs = (1,k)$
for some integer~$k$, then, any of the 
curves~$z^2 = \lambda (z^1)^k$ could be taken to be 
the $z^1$-axis in~$Z_\hs$-linearizing coordinates.
They all have the same tangent space at the origin,
so at most one of them could be geodesic for a given
gradient K\"ahler Ricci soliton~$g$ defined near~$0\in\bbC^2$
with associated holomorphic vector field~$Z_\hs$.  

However, if one uses $g$-normalized coordinates as provided
by Theorem~\ref{thm: normedcoords}, there is a 
canonical~$\bbC^{\mu_i}\subset\bbC^n$ associated 
to the eigenvalue~$h_i$ of multiplicity~$\mu_i$ 
by the equations~$z^j=0$ when~$h_j\not=h_i$.  
It is still not clear to me whether this canonical 
subspace is totally geodesic unless~$h_i$ satisfies 
the `no supermultiples' condition.
\end{remark}

\subsubsection{Growth of $f$ in linearizing coordinates}
\label{sssec: fgrowthinz}
Now that global linearizing coordinates are available, 
it makes sense to ask about the growth of the metric~$g$
and its related quantities in those coordinates. 

One particularly useful quantity to estimate will be
the size of~$|\nabla f|^2(z)$ as $|z|\to\infty$.  
Note that, because of~\eqref{eq: phi2ndpos}, the function
$|\nabla f|^2$ is strictly increasing on the nonconstant 
flow lines of~$\nabla f$.  On the other hand, $|\nabla f|^2
= 2h - R(g)$ is bounded by~$2h$.  Define
\be
\lambda_- = \liminf_{|z|\to\infty} |\nabla f|^2(z) > 0
\qquad\text{and}\qquad
\lambda_+ = \sup_{z} |\nabla f|^2(z) \le 2h.
\ee 

\begin{proposition}\label{prop: fgrowthinz}
For any~$r>0$, there exist constants~$a_1>0$, $a_2>0$, 
$b_1$, and~$b_2$ such that, for all~$z\in\bbC^n$ 
with~$|z|\ge r$,
\be
a_1\,\log|z|+b_1 \le f(z) \le a_2\,\log|z|+b_2\,.
\ee
Explicitly, one can take
\be\label{eq: a1a2def}
a_1 = \frac1{h_n}\inf_{|z|=r}|\nabla f(z)|^2(z) > 0
\qquad\text{and}\qquad
a_2 = \frac{\lambda_+}{h_1} \le \frac{2h}{h_1}.
\ee
\end{proposition}

\begin{proof}
Fix~$r>0$ and note that there exist constants~$m_r>0$ and~$M_r>0$
such that
\be
m_r \le f(z) \le M_r\qquad\text{when $|z|=r$.}
\ee
Moreover, taking~$a_1$ and~$a_2$ as defined in
\eqref{eq: a1a2def} and using the fact that~$|\nabla f|^2(z)$
and $|z|$ both increase along the flow lines of~$\nabla f$, 
one sees that
\be
h_n\,a_1\le |\nabla f(z)|^2\le h_1\,a_2 \qquad\text{when $|z|\ge r$.}
\ee

Now, the flow of~$\nabla f = 2\text{Re}(Z)$ in $Z$-linearizing 
coordinates is
\be
\exp_{t\nabla f}(z^1,\ldots,z^n)
=(\E^{h_1t}z^1,\ldots,\E^{h_nt}z^n),
\ee 
so, since~$0<h_1\le\cdots\le h_n$, it follows that
\be
\E^{h_1t}|z| \le \left|\exp_{t\nabla f}(z^1,\ldots,z^n)\right|
\le \E^{h_nt}|z|.
\ee
In particular, it follows that, for~$t\ge0$.
\be
t\le \frac1{h_1}
\left(\log\left(\left|\exp_{t\nabla f}(z^1,\ldots,z^n)\right|\right) 
- \log|z|\right).
\ee
and
\be
\frac1{h_n}
\left(\log\left(\left|\exp_{t\nabla f}(z^1,\ldots,z^n)\right|\right) 
- \log|z|\right)
\le t.
\ee

On the other hand, since~$\Lie_{\nabla f}f = |\nabla f|^2$,
it follows that
\be
f(z) + h_n\,a_1\,t \le
f\bigl(\exp_{t\nabla f}(z^1,\ldots,z^n)\bigr)
\le f(z) + h_1\,a_2\,t
\ee
for all~$t\ge0$ and~$z$ satisfying~$|z|=r$.  Combining
this with the above inequality gives, for all~$t\ge0$
and~$z$ satisfying~$|z|=r$,
\be
f\bigl(\exp_{t\nabla f}(z^1,\ldots,z^n)\bigr)
- a_2
\log\left|\exp_{t\nabla f}(z^1,\ldots,z^n)\right|
\le f(z) - a_2\,\log|z|.
\ee
Since every~$w\in\bbC^n$ with~$|w|\ge r$ is of the
form~$w = \exp_{t\nabla f}(z)$ for some~$t\ge0$
and $z$ with~$|z|=r$, it follows that
\be
f(w) \le a_2\,\log|w|+\left(M_r-a_2\,\log r\right)
\ee 
for all~$w\in\bbC^n$ with~$|w|\ge r$.  
Thus, taking~$b_2 = M_r-a_2\,\log r$ 
verifies the claimed upper bound on~$f$.

The lower bound follows by combining the lower bound on~$t$
with the lower bound on~$f$:
\be
m_r + a_1
\left(\log\left(\left|\exp_{t\nabla f}(z^1,\ldots,z^n)\right|\right) 
- \log|z|\right) \le f\bigl(\exp_{t\nabla f}(z^1,\ldots,z^n)\bigr),
\ee
which gives
\be
\left(m_r - a_1\log r\right) 
+ a_1\log|w|  \le f\bigl(w\bigr),
\ee
for all~$w\in\bbC^n$ with~$|w|\ge r$.
\end{proof}

Note that, as a function of~$r$, the expression~$a_1$
defined in~\eqref{eq: a1a2def} is
increasing and its limit as~$r\to\infty$ is~$\lambda_-/h_n$.

\begin{corollary}\label{cor: fandlogzbounds}
For any~$\varepsilon>0$, there exists~$r>0$
such that, for~$z\in\bbC^n$ with~$|z|\ge r$,
\be
\left(\frac{\lambda_-}{h_n} - \varepsilon\right)\,\log|z| 
\le f(z) \le 
\left(\frac{\lambda_+}{h_1} + \varepsilon\right)\,\log|z|.
\ee
In particular, there exist constants~$b_1>0$ and~$b_2>0$
such that, for all~$z\in\bbC^n$ with~$|z|\ge r$,
\be
b_1\log|z| \le d(z,p) \le b_2\log|z|.
\ee \qed
\end{corollary}

\begin{proof}
The first statement follows by elementary reasoning
from Proposition~\ref{prop: fgrowthinz} while the 
second follows by combining the first with 
Corollary~\ref{cor: flingrowth}.
\end{proof}

Note that~Corollary~\ref{cor: fandlogzbounds} implies that
the ratio~$f(z)/\log|z|$ is bounded above and has
a positive lower bound as $|z|\to\infty$.  Set
\be
\mu_- = \liminf_{|z|\to\infty}\frac{f(z)}{\log|z|}
\qquad\text{and}\qquad
\mu_+ = \limsup_{|z|\to\infty}\frac{f(z)}{\log|z|}.
\ee
Then Corollary~\ref{cor: fandlogzbounds} implies
\be
\frac{\lambda_-}{h_n}\le\mu_- \le \mu_+ \le \frac{\lambda_+}{h_1}.
\ee

\begin{proposition}\label{prop: foverlogzbounds}
One has the bounds $\mu_- \le 2n \le \mu_+$, in other words
\be
\liminf_{|z|\to\infty}\frac{f(z)}{\log|z|}
\le 2n \le 
\limsup_{|z|\to\infty}\frac{f(z)}{\log|z|}\,.
\ee
\end{proposition}

\begin{proof}
Suppose these bounds do not hold and let~$R>0$ 
be fixed large enough so that there
exist positive constants~$a_1$ and $a_2$
where either~$a_2 < 2n$ or else~$a_1 > 2n$
and positive constants $b_1$ and~$b_2$ so that
\be
a_1\log|z| \le f(z) \le a_2\log|z|
\ee 
and
\be
b_1\log|z| \le d(z,0) \le b_2\log|z|
\ee
hold whenever~$|z|\ge R$.   
(Remember that, in these linearizing coordinates~$p=0$.)  

Let~$M>0$ be sufficiently large that~$d(z,0)\le M$ 
when~$|z|\le R$, and consider any real number~$\rho$ 
that is larger than both~$\log R$ and~$M/b_2$.
 
Consider the $g$-metric ball~$B_{b_1 \rho}(0)$.
Since~$d(z,0)\le b_1\rho$ for~$z\in B_{b_1\rho}(0)$, 
it follows that either~$|z|\le R$ or~$b_1\log|z|\le b_1\rho$,
i.e.,~$|z| \le \E^\rho$.  Since~$\E^\rho > R$, in either case
it follows that~$|z| \le \E^\rho$.  Thus~$B_{b_1\rho}(0)$
is contained in the flat metric ball~$B^0_{\E^\rho}(0)$.

On the other hand, if~$|z|\le \E^\rho$,
then either~$|z|\le R$ or else $d(z,0)\le b_2\rho$.
In the former case, $d(z,0)\le M \le b_2\rho$, by
construction.
In either case, $z$ lies in the $g$-metric 
ball~$B_{b_2 \rho}(0)$.  

Thus, one has inclusions
\be
B_{b_1 \rho}(0) 
\subseteq B^0_{\E^\rho}(0)
\subseteq B_{b_2 \rho}(0). 
\ee

Now, the volume form for~$g$ on~$\bbC^n$ is
\be
\text{vol}_g = \E^{-f}\,\text{vol}_0
\ee
where~$\text{vol}_0=\iC^{n^2}2^{-n}\d z\w \d\zbar$ 
is the volume form of the flat metric on~$\bbC^n$.  

Consequently, the $g$-volume of the $g$-metric ball~$B_{b_2 \rho}(0)$ 
is at least as large as the $g$-volume of the 
flat metric ball~$B^0_{\E^\rho}(0)$ which is given
by the integral
\be
\begin{aligned}
\int_{|z|\le \E^\rho} \E^{-f}\,\text{vol}_0
&= \int_{|z|\le R} \E^{-f}\,\text{vol}_0
 + \int_{|z|= R}^{|z|= \E^\rho} \E^{-f}\,\text{vol}_0\\
&\ge  \int_{|z|\le R} \E^{-f}\,\text{vol}_0
 + \int_{|z|= R}^{|z|= \E^\rho} |z|^{-a_2}\,\text{vol}_0\\
&= \int_{|z|\le R} \E^{-f}\,\text{vol}_0
 + \text{vol}(S^{2n-1})
   \int_{s= R}^{s=\E^\rho} s^{2n-1-a_2}\,\d s\\
\end{aligned}
\ee
Now, if~$a_2<2n$, then the above would imply
\be\label{eq: lowervolbound}
\text{vol}\bigl(B_{b_2 \rho}(0),g\bigr)
\ge 
\int_{|z|\le R} \E^{-f}\,\text{vol}_0
 + \frac{\text{vol}(S^{2n-1})}{2n-a_2}
\bigl(\E^{(2n-a_2)\rho}-R^{2n-a_2}\bigr).
\ee 
However, because~$g$ has positive Ricci curvature,
by the Bishop Comparison Theorem~\cite[Theorem 1.3]{MR97d:53001}
the $g$-volume of~$B_{b_2 \rho}(0)$ is bounded
by a constant times~$\rho^{2n}$. Obviously, such
a bound is not compatible with~\eqref{eq: lowervolbound} 
for all $\rho$ sufficiently large. Thus, $a_2\ge 2n$.

In the other direction, the $g$-volume of the $g$-metric 
ball~$B_{b_1 \rho}(0)$ 
is at most as large as the $g$-volume of the 
flat metric ball~$B^0_{\E^\rho}(0)$, which obeys
the upper bound
\be
\begin{aligned}
\int_{|z|\le \E^\rho} \E^{-f}\,\text{vol}_0
&= \int_{|z|\le R} \E^{-f}\,\text{vol}_0
 + \int_{|z|= R}^{|z|= \E^\rho} \E^{-f}\,\text{vol}_0\\
&\le  \int_{|z|\le R} \E^{-f}\,\text{vol}_0
 + \int_{|z|= R}^{|z|= \E^\rho} |z|^{-a_1}\,\text{vol}_0\\
&= \int_{|z|\le R} \E^{-f}\,\text{vol}_0
 + \text{vol}(S^{2n-1})
   \int_{s= R}^{s=\E^\rho} s^{2n-1-a_1}\,\d s\\
\end{aligned}
\ee
If~$a_1 > 2n$, then this would imply
\be\label{eq: uppervolbound}
\text{vol}\bigl(B_{b_1 \rho}(0),g\bigr)
\le 
\int_{|z|\le R} \E^{-f}\,\text{vol}_0
 + \frac{\text{vol}(S^{2n-1})}{a_1 - 2n}
\bigl(R^{2n-a_1} - \E^{(2n-a_1)\rho}\bigr).
\ee 
and the right hand side is bounded as
a function of~$\rho$.  Thus, 
$\text{vol}\bigl(B_{b_1 \rho}(0),g\bigr)$
would be bounded, independent of~$\rho$, which,
because~$g$ is complete and of positive Ricci
curvature on the noncompact manifold~$\bbC^n$, 
violates Theorem~4.1 of~\cite{MR97d:53001},
which asserts that~$g$ must have at least linear
volume growth.  Thus~$a_1\le 2n$.
\end{proof}

\begin{remark}[Growth of~$f$ in examples]\label{rem: exsoffgrwth}
In the case of Hamilton's soliton (Example~\ref{ex: Hamsoliton})
and, more generally Cao's soliton (Example~\ref{ex: Caosoliton}),
one has~$h_1 = h_n$ and~$\lambda_-=\lambda_+=2nh_1$,
so equality holds in the bounds 
of Proposition~\ref{prop: foverlogzbounds}.  

On the other hand for the product examples (Example~\ref{ex: products}),
\be
f(z) = \sum_{k=1}^n \,\log\bigl(1 + (h_k/c_k)|z^k|^2\bigr)
\ee
which satisfies
\be
\liminf_{|z|\to\infty}\frac{f(z)}{\log|z|} = 2
\qquad\text{while}\qquad 
\limsup_{|z|\to\infty}\frac{f(z)}{\log|z|} = 2n.
\ee
In particular, note that this implies~$\lambda_-\le 2h_n< 2h$.
\end{remark}

\begin{remark}[Growth of the potential~$\phi$]
Let~$\phi$ be the~$\bbT_\hs$-invariant potential for~$g$,
i.e., $\Omega = \frac\iC2\,\p\bar\p\phi$, 
and assume that~$\phi$ is normalized so that~$\phi(0)=0$.

Since~$\Lie_{\nabla f}\phi = f$, it follows 
that~$\phi$ is determined in terms of~$f$ and that 
Corollary~\ref{cor: fandlogzbounds} implies growth bounds 
for~$\phi$ as well.  For example, one sees that there
exist positive constants~$r$, $c_1$, and~$c_2$ so that,
whenever~$|z|\ge r$, one has
\be
c_1\,\bigl(\log|z|\bigr)^2 
 \le \phi(z) \le 
c_2\,\bigl(\log|z|\bigr)^2.
\ee

It should be possible to derive~$C^2$-bounds on~$\phi$
(and hence on~$g$) using the fact that~$\phi$ satisfies
an elliptic Monge-Amp\'ere equation, but I do not see,
at present, a good way to do this so as to get any useful 
information.
\end{remark}

\section{The toric case}\label{sec: toric}

In this last section, some remarks will be made about
the reduction of the gradient K\"ahler Ricci soliton
equation in the `toric' case, which will now be defined.

Throughout this section,~$\bbT^n$ will denote the
maximal abelian subgroup of~$\Un(n)$ that consists 
of diagonal matrices.  Although there is no symplectic
form specified on~$\bbC^n$, 
the mapping~$\mu_n:\bbC^n\to\bbR^n$ defined by
\be
\mu_n(z^1,\ldots,z^n) = \bigl(|z^1|^2,\ldots,|z^n|^2\bigr)
\ee
will sometimes be referred to as the `momentum mapping' 
of~$\bbT^n$.

\begin{definition}[Toric metrics]
A $\bbT^n$-invariant K\"ahler metric~$g$ that is
defined on a connected $\bbT^n$-invariant open 
neigborhood of~$0\in\bbC^n$ will be said to be 
\emph{toric}.
\end{definition}

\begin{remark}[Toric ubiquity]\label{rem: toricubiq}
While, at first glance, the toric condition seems to be 
rather special, note that any gradient K\"ahler Ricci
soliton~$g$ on a neighborhood of~$0\in\bbC^n$ 
that has~$(Z,\Upsilon)=(Z_\hs,\d z)$ as 
its associated holomorphic data is invariant under 
the torus~$\bbT_\hs$.  If~$\hs$ is `generic' in the sense 
that the real numbers~$h_1,\ldots,h_n$ are linearly 
independent over~$\bbQ$, then~$\bbT_\hs = \bbT^n$ and 
hence~$g$ is toric.  

Thus, in some sense, the toric case is `generic'
among complete gradient K\"ahler Ricci solitons
with positive Ricci curvature. 
\end{remark}

\subsection{Symmetry reduction in the toric case}
Assuming an $n$-torus symmetry allows one
to reduce the number of independent variables in the 
gradient K\"ahler Ricci soliton equation~\eqref{eq: MAcond}.

\begin{proposition}
Let $g$ be a toric gradient K\"ahler Ricci solition
defined on a connected open neighborhood of~$0\in\bbC^n$ 
with a nonzero associated holomorphic vector field~$Z$ 
and holomorphic volume form~$\Upsilon$ {\upshape(}defined with 
respect to a Ricci potential~$f$ satisfying~$f(0)=0${\upshape)}. 
Then
\begin{enumerate}
\item The vector field~$Z$ is linearized 
in the coordinates~$z = (z^i)$, so that~$Z = Z_\hs$
for some nonzero~$\hs = (h_1,\ldots, h_n)\in\bbR^n$;
\item The $n$-form~$\Upsilon$ is~$c\,\d z^1\w\cdots\d z^n$
for some nonzero constant~$c$; and
\item $g$ has a unique K\"ahler potential satisfying~$\phi(0)=0$ 
of the form
\be\label{eq: phiasa}
\phi(z) = u\bigl(|z^1|^2,|z^2|^2,\ldots,|z^n|^2\bigr)
\ee
for some real-analytic function~$u$ defined on
an open neighborhood of~$0\in\bbR^n$. 
Moreover, $u$ satisfies the singular 
real Monge-Amp\`ere equation
\be\label{eq: MAeqfora}
\det_{1\le i,j\le n}\left(r^i\frac{\p\hfil}{\p r^i}
\left(r^j\frac{\p u}{\p r^j}\right)\right)\, 
\exp\left(\frac12\sum_{j=1}^n {h_j}\,r^j\frac{\p u}{\p r^j}\right)
= |c|^2\,r^1 r^2\cdots r^n\,.
\ee
where
\be\label{eq: poscondsfora}
\ds\prod_{j=1}^n\frac{\p u}{\p r^j}(0) = |c|^2
\qquad\text{and}\qquad
\ds\frac{\p u}{\p r^j}(0)>0,\qquad 1\le j\le n.
\ee
\end{enumerate}

Conversely, for any nonzero~$\hs\in\bbR^n$ 
and any nonzero complex constant~$c$, 
if~$u$ is a real-analytic function defined on an 
open neighborhood of~$0\in\bbR^n$ 
that satisfies~\eqref{eq: MAeqfora} and~\eqref{eq: poscondsfora}, 
then the function~$\phi$ defined on a $\bbT^n$-invariant neighborhood
of~$0\in\bbC^n$ by~\eqref{eq: phiasa} 
is the K\"ahler potential of a toric 
gradient K\"ahler Ricci soliton on the open neighborhood
of~$0\in\bbC^n$ where it is strictly pseudo-convex.
\end{proposition}

\begin{proof}
To begin with, let me point out a fact that will 
be used several times in the following argument:  
Any $\bbT^n$-invariant holomorphic function defined
on a connected open neighborhood of~$0\in\bbC^n$ is
constant there.  This follows, for example, by
examining the effect of~$\bbT^n$ on the individual
terms in the power series of such a function.
 
Now, since~$g$ is toric, its associated holomorphic
vector field~$Z$ is invariant under
the action of~$\bbT^n$ and hence must vanish at~$0\in\bbC^n$ 
and commute with each of the scaling vector fields~$Z_i 
= z^i\frac{\p\hfill}{\p z^i}$.  It follows easily
that~$Z=Z_\hs$ for some~$\hs\in\bbR^n$.  (For the definition
of~$Z_\hs$, see~\eqref{eq: Z_hsdef}.)

Let~$f$ be the unique $\bbT^n$-invariant Ricci potential
for~$g$ that satisfies~$f(0)=0$ and let~$\Upsilon$ 
be a holomorphic volume form associated to~$g$ and~$f$.
Since~$\Upsilon$ is uniquely determined up to a
complex number of modulus~$1$, it follows that, under
the action of~$\bbT^n$, $\Upsilon$ must transform
multiplicitively by a character of~$\bbT^n$.  It then
follows easily that~$\Upsilon= c\,\d z$ 
for some nonzero constant~$c$.

Let~$\phi$ be the unique $\bbT^n$-invariant K\"ahler
potential for~$g$ that satisfies~$\phi(0)=0$.  As has
already been remarked, $\phi$ is real-analytic and so
can be expanded as a convergent power series in the 
variables~$z^i$ and~$\zbar^i$.  However, $\bbT^n$-invariance
evidently implies that this power series can be collected
in terms of the quantities~$r^i = |z^i|^2$.  Thus,
the existence of a function~$u$ satisfying~\eqref{eq: phiasa}
follows.

As argued in~\S\ref{ssec: nonsingextprobs}, 
the quantity~$2f - \p\phi(Z_\hs)$ 
is a holomorphic function on a neighborhood
of~$0\in\bbC^n$.  By construction, it, too, 
is $\bbT^n$-invariant and vanishes at~$0\in\bbC^n$,
which implies that it vanishes identically. 
Thus, $\p\phi(Z_\hs) = \d\phi(X_\hs) = 2f$.

The rest of the argument follows by substituting the
formula~\eqref{eq: phiasa} into~\eqref{eq: MAcond},
multiplying by~$r^1\cdots r^n$, 
and rearranging terms, 
which gives~\eqref{eq: MAeqfora}.  

Note that the stated positivity conditions on the 
first derivatives of~$u$ are needed
in order that the corresponding~$\phi$ be strictly 
pseudo-convex in a neighborhood of~$0\in\bbC^n$ and
the relation with~$|c|^2$ follows by computing the
coefficient of~$r^1\cdots r^n$ in the power series
expansion of the left hand side of~\eqref{eq: MAeqfora}.

The converse statement follows by computation.
\end{proof}

\begin{remark}[Normalizations]
Given a solution~$u$ to~\eqref{eq: MAeqfora} that
satisfies~$u(0)=0$, one can obviously scale in the 
individual coordinates so as to arrange that
\be
\phi = r^1 + \cdots + r^n + O\bigl(|r|^2\bigr),
\ee
thereby reducing to the case~$|c|=1$, so it suffices
to consider this case.  Note also that the
resulting K\"ahler soliton~$g$ is already 
in the normalized form guaranteed 
by~Theorem~\ref{thm: normedcoords}.
\end{remark}

\begin{remark}[pseudo-convexity of toric potentials]
A $\bbT^n$-invariant function~$\phi$ 
of the form~\eqref{eq: phiasa}, i.e., $\phi = u\circ\mu_n$
for some~$u$ defined on a domain~$V\subset\bbR^n$,
is strictly pseudo-convex on the domain~$(\mu_n)^{-1}(V)\subset\bbC^n$ 
if and only if the symmetric matrix
\be
\begin{pmatrix} \ds
\delta_{ij}\,\frac{\p u}{\p r^j}
+ \sqrt{r^ir^j}\,\frac{\p^2 u}{\p r^i\p r^j}
\end{pmatrix}
\ee
is positive definite on the part of~$V$ that lies
in the orthant defined by the inequalities~$r^i\ge0$.
\end{remark}

\subsubsection{A singular initial value problem}
Although~\eqref{eq: MAeqfora}
is singular along the hypersurfaces~$r^i=0$ 
in~$\bbR^n$, it turns out that the methods of
G\'erard and Tahara \cite{MR2001c:35056} can be used
to prove an extension theorem.

\begin{theorem}
Let~$v$ be a real-analytic function 
on an open subset~$V\subset\bbR^{n-1}$ 
with the property that~$\psi=v\circ\mu_{n-1}$ 
is strictly pseudo-convex 
on~$(\mu_{n-1})^{-1}(V)\subset\bbC^{n-1}$.

Then there exists an open neighborhood~$U\subset\bbR^n$
of~$V\times\{0\}$ and a real-analytic function~$u$ on~$U$
with the properties
\begin{enumerate}
\item $u(r^1,\ldots,r^{n-1},0) = v(r^1,\ldots,r^{n-1})$
  for~$(r^1,\ldots,r^{n-1})\in V$;
\item $u$ satisfies~\eqref{eq: MAeqfora} with~$|c|=1$; and
\item $\phi = u\circ\mu_n$ is strictly pseudo-convex 
  on~${\mu_n}^{-1}(U)\subset\bbC^n$.
\end{enumerate}

Moreover,~$u$ is locally unique in the sense that any
for any other pair~$(\tilde U, \tilde u)$ with these 
properties, there is an open neigborhood~$W$ of~$V\times\{0\}$
contained in~$U\cap\tilde U$ such that~$u$ and~$\tilde u$ 
agree on~$W$.
\end{theorem}

\begin{proof}
For the sake of clarity, write~$t = r^n$ and let
the lower case latin indices run from~$1$ to~$n{-}1$.  Then 
after dividing both sides of~\eqref{eq: MAeqfora} (with~$|c|=1$)
by~$r^1\cdots r^{n-1}$ and the exponential factor, 
this equation takes the form
\be\label{eq: tsingMAforu1}
\det\left(
\begin{matrix}
\ds \delta_{ij}\frac{\p u}{\p r^i} + r^j\,\frac{\p^2 u}{\p r^i\p r^j}
& \ds \frac{\p (tu_t)}{\p r^i}\\[10pt]
\ds r^j\frac{\p (tu_t)}{\p r^j} & (t\,\p_t)^2 u
\end{matrix}\right)\, 
= t\,\E^{\ds\left(-\frac{h_n}2\,(tu_t)
-\frac12\sum_{j=1}^{n-1} {h_j}\,r^j\frac{\p u}{\p r^j}\right)}.
\ee
Note the first crucial aspect of this equation, which is
that the $t$-derivatives of~$u$ occur as either~$tu_t$ or
$t(tu_t)_t = (t\p_t)^2 u$, i.e., as the `regular singular'
versions of the $t$-derivatives at~$t=0$.  

Expanding the left hand side of~\eqref{eq: tsingMAforu1}
along the last column shows that this equation can be
written in the form
\be\label{eq: tsingMAforu2}
\begin{aligned}
\det\left(\delta_{ij}\frac{\p u}{\p r^i} 
+ r^j\,\frac{\p^2 u}{\p r^i\p r^j}\right)\,\bigl((t\,\p_t)^2 u\bigr)
&= t\,\E^{\ds\left(-\frac{h_n}2\,(tu_t)
-\frac12\sum_{j=1}^{n-1} {h_j}\,r^j\frac{\p u}{\p r^j}\right)}
\\
&\qquad+Q_{ij}\left(r,\frac{\p u}{\p r},\frac{\p^2 u}{\p r^2}\right)
\,\frac{\p (tu_t)}{\p r^i}\,\frac{\p (tu_t)}{\p r^j}
\end{aligned}
\ee
where~$Q_{ij}=Q_{ji}$ are certain polynomials in the variables~$r^i$
and the first and second derivatives of~$u$ with respect to
the variables~$r^i$. 

In particular, note that the right hand side
of~\eqref{eq: tsingMAforu2} is an entire analytic function
of the variables~$r^i$ and~$t$, the first and second derivatives 
of~$u$ with respect to the variables~$r^i$, the expression~$tu_t$
and its first derivatives with respect to the~$r^i$.  

\emph{In what
follows, it will be particularly important that this right
hand side is also in the ideal generated by~$t$ and the quadratic 
expressions~$\frac{\p (tu_t)}{\p r^i}\,\frac{\p (tu_t)}{\p r^j}$.}

Now, set~
\be\label{eq: uisvplusz}
u(r^1,\ldots,r^{n-1},t) 
= v(r^1,\ldots,r^{n-1})
+z(r^1,\ldots,r^{n-1},t)
\ee
and define
\be
F_{ij}(r^1,\ldots,r^{n-1},t) 
= \delta_{ij}\frac{\p v}{\p r^i} + r^j\,\frac{\p^2 v}{\p r^i\p r^j}
\ee
Note that, by hypothesis, $\det\bigl(F_{ij}(r,0)\bigr)\not=0$
for~$r\in V\subset\bbR^{n-1}$.  In particular, the expression
\be
\det\left(F_{ij}(r,t) +
\delta_{ij}\frac{\p z}{\p r^i} 
+ r^j\,\frac{\p^2 z}{\p r^i\p r^j}\right),
\ee
which is what the coefficient of~$(t\p_t)^2u$
on the left hand side of~\eqref{eq: tsingMAforu2} 
becomes when one substitutes~$u = v + z$ into
that equation, is an analytic expression in~$r\in V$, 
$t$, and the partials of~$z$ that is non-vanishing on~$V$
when one sets~$t=z=0$.

Thus, substituting~$u=v+z$ into~\eqref{eq: tsingMAforu2}
and dividing by the determinant factor yields 
an equation for~$z$ of the form
\be\label{eq: tsingMAforz}
(t\,\p_t)^2 z
=E\left(r,t,z,\frac{\p z}{\p r^i},tz_t,
      \frac{\p^2 z}{\p r^i\p r^j},\frac{\p (tz_t)}{\p r^i}\right)
\ee
where the function~$E$ is
\begin{enumerate}
\item real-analytic 
on an open neighborhood of~$V\times\{0\}$
in~$V\times\bbR\times\bbR^{1+n+\frac12n(n+1)}$ and
\item in the ideal generated by~$t$ 
and the products of \emph{pairs} of the last~$(n{-}1)$ variables
(i.e., the `slots' containing the entries~$\frac{\p (tz_t)}{\p r^i}$).
\end{enumerate} 

Now, turning to Chapter~8 of G\`erard and Tahara~\cite{MR2001c:35056}, 
one sees that~\eqref{eq: tsingMAforz}
is of the form to which their~Theorem~8.0.3 applies.%
\footnote{While I do not want to state their full theorem 
here, I will give the gist:  The two properties listed 
for the function~$E$ are easily seen to imply that there
exists a unique formal power series solution
of the form~$z(r,t)=z_1(r)\,t + z_2(r)\,t^2 + \cdots$ 
to~\eqref{eq: tsingMAforz}.  The main import of
the quoted Theorem~8.0.3 is that this series actually
converges to an analytic solution on some open neighborhood
of~$V\times\{0\}$.  (The need for a theorem is caused by
the singularity at~$t=0$, which renders the standard 
method of majorants ineffective in proving the convergence 
of the formal series.)} 
Consequently, \eqref{eq: tsingMAforz} has a unique 
real-analytic solution~$z(r,t)$ (defined on some neighborhood 
of~$V\times\{0\}\subset\bbR^n$) that satisfies the initial condition
\be
z(r^1,\ldots,r^{n-1},0) = 0\qquad\text{for~$(r^1,\ldots,r^{n-1})\in V$.}
\ee
Using this solution~$z$ to define~$u$ via~\eqref{eq: uisvplusz},
one sees that~\eqref{eq: tsingMAforu1} has a correspondingly unique 
real-analytic solution satisfying the initial condition
\be
u(r^1,\ldots,r^{n-1},0) = v(r^1,\ldots,r^{n-1})
\qquad\text{for~$(r^1,\ldots,r^{n-1})\in V$,}
\ee
as claimed.  The existence of an open neighborhood~$U$
of~$V\times\{0\}$ such that~$\phi=u\circ\mu_n$
is strictly pseudo-convex on~$(\mu_n)^{-1}(U)\subset\bbC^n$
is routine.
\end{proof}

\begin{corollary}[Singular initial value problem for toric solitons]
\label{cor: singIVP}
Let~$g'$ be a real-analytic toric K\"ahler metric 
on a $\bbT^{n-1}$-invariant, connected 
open neighborhood~$V\subset\bbC^{n-1}$ of~$0$.  

Then, for any~$\hs\in\bbR^n$ there exists a $\bbT^n$-invariant 
open neighborhood~$U_\hs\subset\bbC^n$
of~$V\times\{0\}$ and a toric gradient K\"ahler Ricci 
soliton~$g_\hs$ on~$U_\hs$ whose pullback to~$V$ is~$g'$,
whose associated vector field is~$Z_\hs$, and whose
associated holomorphic volume form with respect
to its~$\bbT^n$-invariant Ricci potential~$f_\hs$
vanishing at~$0\in\bbC^n$  
is~$\Upsilon = \d z^1\w\cdots\w \d z^n$.

Moreover, $g_\hs$ is locally unique in that any 
extension of~$g'$ with these properties agrees with~$g_\hs$
on some open neighborhood of~$V\times\{0\}$. \qed
\end{corollary}

\begin{remark}[Contrast in initial value problems]
Note that Corollary~\ref{cor: singIVP} has a very 
different character from~Theorem~\ref{thm: locgen}.
Not only is the nature of the initial data different, 
but, in the case of Corollary~\ref{cor: singIVP}, 
one is imposing initial conditions along a submanifold
that is everywhere \emph{tangent} to the holomorphic vector
field~$Z = Z_\hs$, rather than everywhere \emph{transverse}.
The difference, of course, is that Corollary~\ref{cor: singIVP}
addresses a singular initial value \textsc{pde} problem 
that is, in many ways the analogue of the sort of \textsc{ode}
problem one encounters in the theory of regular singular
points of \textsc{ode}.  

Because the generalization of the \textsc{ode} concept 
of `regular singular point' to the case of~\textsc{pde}
is very delicate (cf.~the book of G\`erard and Tahara),
it is somewhat remarkable that this theory actually applies 
in this case.
\end{remark}

\subsubsection{A Lagrangian formulation}
While the reduced equation~\eqref{eq: MAeqfora} 
is singular along the hypersurfaces~$r^i=0$, 
it is regular on the open simplicial cone
defined by~$r^i>0$.  
Indeed, setting~$r^i = \E^{\rho^i}$, 
the equation~\eqref{eq: MAeqfora} with~$|c|^2=1$ 
can be written in the form
\be\label{eq: MAasurho}
\det_{1\le i,j\le n}\left(\frac{\p^2 u}{\p\rho^i\p\rho^j}\right)\, 
\E^{\left(\frac{h_1}2\,\frac{\p u}{\p\rho^1}
+ \cdots + \frac{h_n}2\,\frac{\p u}{\p\rho^n}\right)}
= \E^{\rho^1+\cdots+\rho^n}\,.
\ee
Setting~$u_i = \frac{\p u}{\p\rho^i}$, this can be further
rewritten into the form
\be
\E^{\left(\frac{h_1}2\,u_1+\cdots+\frac{h_n}2\,u_n\right)}
\,\d u_1\w \cdots\w \d u_n
= \E^{\rho^1+\cdots+\rho^n}\,\d \rho^1\w\cdots\w \d\rho^n.
\ee
Thus, on~$\bbR^{2n+1}$ with coordinates~$u,\rho^i,u_i$, if
one defines the contact form
\be
\theta = \d u - u_i\,\d\rho^i
\ee
and the closed $\theta$-primitive%
\footnote{If~$(M^{2n+1},\theta)$ is a contact manifold
of dimension~$2n{+}1$, then an $n$-form~$\Psi$ on~$M$
is said to be \emph{$\theta$-primitive} if~$\d\theta\w\Psi
\equiv 0 \mod\theta$.}
 $n$-form
\be
\Psi 
= \E^{\left(\frac{h_1}2\,u_1+\cdots+\frac{h_n}2\,u_n\right)}
      \,\d u_1\w \cdots\w \d u_n
     - \E^{\rho^1+\cdots+\rho^n}\,\d \rho^1\w\cdots\w \d\rho^n,
\ee
Then the solutions of the original equation~\eqref{eq: MAeqfora}
correspond to the integral manifolds of the Monge-Amp\`ere ideal
\be
\cI = \la \theta,\d\theta,\Psi\ra.
\ee

Since~$\Psi$ is closed and~$\d\theta\w\Psi=0$, 
the $(n{+}1)$-form~$\Pi = \theta\w\Psi$
is closed and hence is the Poincar\'e-Cartan form
(see~\cite{MR04g:58001}) of a contact Lagrangian 
for the function~$u$. In particular, it follows by
Noether's Theorem that the symmetries of the 
Poincar\'e-Cartan form give conservation laws 
for the reduced equation.   

This is interesting because this system turns out to
have a number of symmetries that are not apparent from
the symmetries of the original equation.

\begin{remark}[Affine symmetries and equivalences]
For example, consider the affine transformations on~$\bbR^{2n+1}$
of the form
\be
\begin{aligned}
\bar u &= s\,u + a_i B^i_k\,\rho^k + c\,,\\
\bar u_i &= A_i^j\,u_j + a_i\,,\\
\bar\rho^i &= B^i_j\,\rho^j + b^i
\end{aligned}
\ee
where~$A^i_j$, $B^i_j$, $s\not=0$, $a_i$, $b^i$, and~$c$ are
real constants satisfying the $n^2 + 2n + 1$ equations
\be
\begin{aligned}
A^j_iB^i_k &= s\,\delta^j_i\,,\\
\sum_i h_i A^j_i & = h_j\quad\text{for $1\le j\le n$}\,,\\
\sum_i B^i_j&=1\quad\text{for $1\le j\le n$}\,,\\
\E^{\left(\frac{h_1}2\,a_1+\cdots+\frac{h_n}2\,a_n\right)}
     \,\det(A) &=  \E^{b^1+\cdots+b^n}\,\det(B)\,.
\end{aligned}
\ee
Such transformations, which constitute a Lie group of 
dimension~$n^2+1$, preserve the forms~$\theta$ and~$\Upsilon$
up to constant multiples and hence preserve the system~$\cI$.

Obviously, the system depends on the vector~$\hs = (h_1,\ldots,h_n)$.
However, by leaving off the second of the above four conditions,
one finds transformations that define equivalences between
any two systems with~$h = h_1+\cdots+h_n\not=0$ and any two 
systems with~$h = h_1+\cdots+h_n=0$ but~$\hs\not=0$.
(The system corresponding to~$\hs=0$ is, of course, the
system that gives Ricci-flat toric K\"ahler metrics.)
\end{remark}

\begin{remark}[Algebraic coordinates]
The function~$u$ is, in some sense, not that important, since
only the derivatives of~$u$ appear in the formula for the
metric.  Thus, one can actually formulate the essential 
part of the exterior differential system as a system on~$\bbR^{2n}$.

Assuming that none of the~$h_i$ are zero, one can coordinatize
the system algebraically as follows:  Set~$v_i = \E^{\frac12 h_i u_i}$.
Then the form~$\Upsilon$ can, after multiplying by a constant
be written in the form
\be
\Upsilon = \d v_1 \w \cdots \w \d v_n 
- \frac{h_1\cdots h_n}{2^n}\,\d r^1 \w \cdots \w \d r^n,
\ee
and the contact condition that~$\d u - u_i\,\d\rho^i = 0$ 
can be replaced by the condition
\be
\sum_{i=1}^n \frac2{h_i} \frac{\d v_i}{v_i} \w \frac{\d r^i}{r^i} = 0.
\ee
\end{remark}

\bibliographystyle{hamsplain}

\providecommand{\bysame}{\leavevmode\hbox to3em{\hrulefill}\thinspace}

\end{document}